\documentclass[11pt]{amsart}

\usepackage[T1]{fontenc}
\usepackage[utf8]{inputenc}
\usepackage{lmodern}
\usepackage{microtype}
\usepackage[a4paper,margin=2.7cm]{geometry}
\usepackage{amsmath,amssymb,amsthm,mathtools}
\usepackage{booktabs,array}
\usepackage{enumitem}
\usepackage{graphicx}
\usepackage[hidelinks]{hyperref}
\hypersetup{pdftitle={A tower of reducible sets and new lower bounds for three-term progression-free sets in F\_p\^{}n and in the integers}}

\newcommand{\Z}{\mathbb{Z}}
\newcommand{\Zm}{\mathbb{Z}/m\mathbb{Z}}
\newcommand{\hf}[1]{\langle #1\rangle}          
\DeclareMathOperator{\vol}{vol}
\newcommand{\file}[1]{\texttt{#1}}

\theoremstyle{plain}
\newtheorem{theorem}{Theorem}[section]
\newtheorem{proposition}[theorem]{Proposition}
\newtheorem{lemma}[theorem]{Lemma}
\newtheorem{corollary}[theorem]{Corollary}

\theoremstyle{definition}
\newtheorem{definition}[theorem]{Definition}
\theoremstyle{remark}
\newtheorem{remark}[theorem]{Remark}
\numberwithin{equation}{section}

\setlist[itemize]{leftmargin=2.2em}
\setlist[enumerate]{leftmargin=2.2em}

\title[A tower of reducible sets]{A tower of reducible sets and new lower bounds for three-term
progression-free sets in $\mathbb{F}_p^n$ and in the integers}
\author{Christian Elsholtz}
\address{Institut f\"ur Analysis und Zahlentheorie, Graz University of Technology, Kopernikusgasse 24,
8010 Graz, Austria}
\author{Arnd Roth}
\address{Wolfson Institute for Biomedical Research, University College London, Gower Street, London WC1E 6BT,
United Kingdom}
\date{29 November 2026}

\begin{document}

\begin{abstract}
We develop an inductive construction of non-mid-point reducible sets in tori of arbitrary dimension. It is a
new type of lifting to higher dimensions, which can improve upon the density obtained by Cartesian products. Starting from a one-dimensional
interval, the construction produces a tower of sets whose slice measures satisfy a linear recursion. We prove
an additional exponential growth compared with the previous constructions. By the lifting
principle of Elsholtz, Proske and Sauermann this yields three-term progression-free subsets of $(\Zm)^n$ of
size $(cm)^{n-O(\log n)}$ with $c=0.554993$, for every fixed $m\ge2$, where $c=\sqrt{7/24}\approx0.5401$ was
known. A quantitative version of the induction preserves the potentials that Elsholtz, Hunter, Proske and
Sauermann use in Behrend's construction, up to arbitrarily small losses of volume. Consequently
$r_3(N)\ge N\cdot2^{-(C+o(1))\sqrt{\log_2N}}$ with $C<2.60685$, improving the previous constant
$2\sqrt{\log_2(24/7)}\approx2.6665$. A different removal order yields larger reducible sets in dimensions two and three, and improves the
construction in dimension four. We also prove upper bounds for reducible sets and for the inductive
construction.
\end{abstract}

\maketitle

\section{Introduction and main results}\label{sec:intro}

A subset $A$ of an abelian group $G$ is called \emph{three-term progression free}, or
\emph{progression-free} for short, if it contains no three distinct elements $x,y,z$ with $x+z=2y$. We write
$r_3(G)$ for the largest size of a progression-free subset of $G$, and $r_3(N)$ for the largest size of a
progression-free subset of $\{1,\dots,N\}$. This paper is about lower bounds for $r_3((\Zm)^n)$, where
$m\ge 2$ is a fixed integer and $n\to\infty$, and about lower bounds for $r_3(N)$.

The upper bounds for the first problem were transformed by the polynomial method. Croot, Lev and
Pach~\cite{CLP} showed that progression-free subsets of $(\Z/4\Z)^n$ are exponentially small, and
Ellenberg and Gijswijt~\cite{EG} proved $r_3(\mathbb{F}_p^n)\le (J_p\,p)^n$ for odd primes $p$, with constants
$J_p<1$ that decrease to $J_\infty=0.84143\ldots$ as $p\to\infty$. For the integers, Kelley and
Meka~\cite{KM} proved $r_3(N)\le N\exp(-c(\log N)^{1/12})$ for some $c>0$, and Bloom and Sisask~\cite{BS}
improved the exponent $\frac1{12}$ to $\frac19$. Recently Raghavan~\cite{Raghavan} improved it further, to
$\frac16$ up to a factor $\log\log N$ in the exponent.

For lower bounds and an odd prime $p$, a classical construction takes the vectors in
$\{0,1,\dots,\frac{p-1}{2}\}^n$ with a prescribed sum of squares. Such a set is progression-free, because
$x+z=2y$ in $\mathbb{F}_p^n$ then holds in $\Z^n$, and the sum of squares is strictly convex. This gives
$r_3(\mathbb{F}_p^n)\ge \bigl(\frac{p+1}{2}\bigr)^{n-O(\log n)}$. The older construction of
Salem and Spencer~\cite{SalemSpencer}, in which every digit of $\{0,1,\dots,\frac{p-1}{2}\}$ occurs equally
often, gives the same base. For small values of a prime $p$ or of a composite $m$ much better constructions are known, see for
example~\cite{ElsholtzPach} for $m=4$ and $m=5$. When $p=3$ this is
the well-known case of caps, see for example~\cite{FunSearch,Tyrrell}. When $p=5$, the first author and
students proved, in work that is not yet published, that one can lift a set of $36$ reducible points from
dimension~$3$ to much higher dimensions, and together with a refinement one can achieve
$r_3(\mathbb{F}_5^n)\ge3.3074^n$ and $r_3((\Z/15\Z)^n)\ge8.785^n$ for large values of $n$. In
the integers the same idea is Behrend's construction~\cite{Behrend}, which gives
\[
  r_3(N)\;\ge\;N\cdot2^{-(C+o(1))\sqrt{\log_2N}}
\]
with $C=2\sqrt2\approx2.8284$. For almost eighty years only the $o(1)$ term in this bound was
improved~\cite{Elkin,GreenWolf}.

Elsholtz, Proske and Sauermann~\cite{EPS} improved the base $\frac{p+1}2$, and Elsholtz, Hunter, Proske and
Sauermann~\cite{EHPS} improved the constant $C$ for the first time. The method of~\cite{EPS} has two parts.
The first is a lifting principle: if a subset $S$ of $(\Zm)^d$ is \emph{non-mid-point reducible}
(Definition~\ref{def:reducible}), which means that its points can be removed one at a time so that no removed
point is the mid-point of two points that are still present, then $(\Zm)^n$, in a much larger dimension $n$,
contains a progression-free subset of size $|S|^{n/d-O(\log n)}$ (Proposition~\ref{prop:eps}). The resulting
constant is $c=(|S|/m^d)^{1/d}$, the density of $S$ per dimension. Intervals of length about $m/2$ in $\Zm$ are
reducible, which recovers the constant $\frac12$. The second part is a construction in dimension~$2$: a
reducible subset of the torus $[0,1)^2=\mathbb{R}^2/\Z^2$ of area $\frac7{24}$, whose discretisations are
reducible subsets of $(\Zm)^2$ of density at least $\frac7{24}$ for every $m$. This gives
$c=\sqrt{7/24}\approx0.5401$. In~\cite{EHPS} it is shown that subsets of the torus $[0,1)^d$ that carry a
suitable function (a \emph{potential}, Definition~\ref{def:potential}) can be used as generalised digits in
Behrend's construction; a set of volume $v$ gives $C=2\sqrt{\log_2(1/c^2)}$ with $c=v^{1/d}$, and the
two-dimensional set gives $C=2\sqrt{\log_2(24/7)}\approx2.6665$.

In this paper we come back to the concept of reducible sets and refine the constructions. In
dimension~$2$ we give a small improvement, a reducible set of area about $0.2943$ instead of
$\frac7{24}\approx0.2917$ (Theorem~\ref{thm:dim3}). Our main result is an induction step
(Theorem~\ref{thm:step}) that turns a suitable reducible set in dimension~$d$ and an arbitrary reducible set
in dimension $d-1$ into a reducible set in dimension $d+1$. It shows that one achieves slowly increasing
improvements as $d\to\infty$: reducible sets $T_d\subseteq[0,1)^d$, one in every dimension, whose constants
$(\vol T_d)^{1/d}$ increase on all floors that we computed and converge to a limit, which we enclose in an
interval of length $10^{-8}$ (Theorem~\ref{thm:main}). This appears to be the first improvement in a
dimension $d>2$; it was not known before whether reducible sets in higher dimensions can give a larger
constant than products of two-dimensional ones. When we apply the construction to Behrend's construction, we
obtain a further improvement of the constant $C$ (Theorem~\ref{thm:behrend}).

\subsection*{A new lifting to higher dimensions}
The lifting principle of~\cite{EPS} passes from a fixed dimension $d$ to all large dimensions, and it
preserves the density per dimension. So the constant $c$ is decided by the reducible sets that are available
in fixed dimensions. The obvious way to obtain reducible sets in higher dimensions is the Cartesian product:
products of reducible sets are reducible, but the density per dimension of a product is a weighted geometric
mean of those of its factors, so a Cartesian product, or a Cartesian power of one set, never improves the constant.

A main point of this paper is a lifting of a different type, from dimension $d$ to dimension $d+1$, which
can improve the density. We call the last coordinate of a point of $[0,1)^d$ its \emph{level}. A
\emph{block} is a reducible subset $B$ of $[0,1)^d$ whose levels are all below $\frac34$ and in which no
three points $x,y,z$ with $x+z\equiv2y$ have levels with $t_x+t_z=2t_y+1$ (Definition~\ref{def:block}). From a block $B$ in
dimension $d$, an arbitrary reducible set $F$ in dimension $d-1$ and a parameter $\beta$ we construct a set
$B^+$ in dimension $d+1$ (Definition~\ref{def:step} and Figure~\ref{fig:step}), and we prove that $B^+$ is
again a block (Theorem~\ref{thm:step}). A product $B\times[0,\frac12)$ uses the same set $B$ at half of the new
levels. The set $B^+$ uses three quarters of the new levels: half of them carry the part of $B$ above a cut
that rises with the level, and a quarter of them carry low slices of $B$ and copies of $F$, each of them
twice. For the blocks of this paper this is more than the product, and since $B^+$ is a block, the step can be
iterated. The interval $[0,\frac12)$ is a block, the first step essentially reproduces the two-dimensional set
of~\cite{EPS} (Remarks~\ref{rem:EPS} and~\ref{rem:floor1}), and the further steps give a tower of reducible
sets $T_2,T_3,T_4,\dots$, one in every dimension, whose density per dimension increases on every floor that we
computed.

The two-dimensional problem has also been attacked with computer searches. In November 2025, in a comment
on a blog post of Tao~\cite{TaoBlog} about the paper~\cite{GGTW}, A.~Z.~Wagner reported that it was the
first problem he tried when he started to work with AlphaEvolve, and that they found many constructions
with area $\sim\frac7{24}$ but nothing better.

\subsection*{Results}

The volume of $B^+$ depends only on the measures of the slices of $B$ and on the volume of $F$, so the volumes
of the floors satisfy a linear recursion, given by a positive compact operator.

\begin{theorem}[The tower]\label{thm:main}\label{thm:limit}
There are non-mid-point reducible sets $T_d\subseteq[0,1)^d$, $d\ge2$ (Definition~\ref{def:tower}), with the
following properties.
\begin{enumerate}
\item[(a)] The ratios $\vol T_{d+1}/\vol T_d$ and the constants $(\vol T_d)^{1/d}$ converge to the same limit
$r_\infty$, and $0.55499308\le r_\infty\le0.55499309$.
\item[(b)] Let $c=0.554993$. For every fixed integer $m\ge 2$ and $n\to\infty$ there is a progression-free
subset $A$ of $(\Zm)^n$ with $|A|\ge(c\,m)^{\,n-O(\log n)}$. In particular
$r_3(\mathbb{F}_p^n)\ge (c\,p)^{\,n-O(\log n)}$ for every fixed prime $p$.
\item[(c)] $\vol T_{1000}>0.55497^{1000}$, so the explicit set $T_{1000}$ gives (b) with $c=0.55497$.
\end{enumerate}
\end{theorem}

Compared with the constant $\sqrt{7/24}$ of~\cite{EPS}, part~(b) is an additional exponential growth: the base
of the lower bound increases by the factor $0.554993/\sqrt{7/24}>1.0276$. The convergence follows from the
theorem of Krein and Rutman. The enclosure of $r_\infty$ and the volume of
$T_{1000}$ rest on computations in exact integer arithmetic (Section~\ref{sec:tower}).

The induction step has a quantitative version for potentials (Theorem~\ref{thm:steppot}): if $B$ and $F$ carry
potentials, then so does $B^+$, after removing parts of arbitrarily small volume. So the floors of the tower
can be used in Behrend's construction.

\begin{theorem}[Integers]\label{thm:behrend}
We have
\[
  r_3(N)\;\ge\;N\cdot 2^{-(C+o(1))\sqrt{\log_2 N}}\qquad\text{with}\qquad
  C=2\sqrt{\log_2\bigl(r_\infty^{-2}\bigr)}<2.60685 ,
\]
where $r_\infty$ is the limit of Theorem~\ref{thm:main}. More generally, every set in
Table~\ref{tab:constants}, of dimension $d$ and volume $v$, gives this bound with
$C=2\sqrt{(2/d)\log_2(1/v)}$, the value in the last column of the table.
\end{theorem}

In the convention with natural logarithms, $r_3(N)\ge N\exp\bigl(-(C'+o(1))\sqrt{\ln N}\bigr)$, the constant is
$C'=C\sqrt{\ln2}<2.17035$, where Behrend has $2\sqrt{2\ln 2}\approx2.3548$ and~\cite{EHPS} has
$\approx2.2200$.

The induction does not exhaust the possibilities in low dimensions. The first two floors leave a region free
below their cuts, and part of it, a notch, can be added if the removal order is changed: the points of the low
levels are no longer removed level by level, but along a family of lines that turn about a fixed point (the
fan order of Definition~\ref{def:fan}).

\begin{theorem}[Low dimensions]\label{thm:dim3}\label{thm:dim4}
There are non-mid-point reducible sets
\begin{enumerate}
\item[(a)] $R_{1/10}\subseteq[0,1)^2$ of area $\frac{209}{800}+\frac9{160}\ln\frac95\approx0.294313$, which is
more than $\frac7{24}\approx0.291667$;
\item[(b)] $T(\frac7{50},\frac25)\cup M^*\subseteq[0,1)^3$ of volume
$\frac{4501789}{22400000}-\frac{4399}{80000}\ln\frac{53}{28}\approx0.165886$;
\item[(c)] $\widetilde T_4\subseteq[0,1)^4$ of volume
$\frac{218969684867}{2400000000000}+\frac{168921}{160000000}\ln\frac95\approx0.0918579$.
\end{enumerate}
They give the conclusion of Theorem~\ref{thm:main}(b) with $c>0.54250$, $c>0.54946$ and $c>0.55052$, and the
bound of Theorem~\ref{thm:behrend} with $C<2.65675$, $C<2.62894$ and $C<2.62468$, respectively.
\end{theorem}

The set in (a) shows that the set of~\cite{EPS} is not the largest reducible subset of the two-dimensional
torus. The sets in (a) and (b) are not blocks, so they cannot serve as floors of the tower, but they can serve
as free sets $F$: the set in (c) is a fourth floor with the free set $R_{1/10}$. Table~\ref{tab:constants}
lists the constants.

\begin{table}[ht]
\centering
\small
\renewcommand{\arraystretch}{1.25}
\begin{tabular}{@{}cllll@{}}
\toprule
$d$ & Set & Volume $v$ & $c=v^{1/d}$ & $C$\\
\midrule
$1$ & $[0,\frac12)$ & $1/2$ & $0.5$ & $2.82843$\\
$2$ & $E(5/12)$, cf.~\cite{EPS} & $7/24\approx 0.29167$ & $0.54006$ & $2.66654$\\
$2$ & $R_{1/10}$ & $\approx 0.29431$ & $0.54250$ & $2.65675$\\
$3$ & $T_3=T(\frac18,\frac25)$ & $15803/96000\approx 0.16461$ & $0.54805$ & $2.63456$\\
$3$ & $T(\frac7{50},\frac25)\cup M^*$ & $\approx 0.16589$ & $0.54946$ & $2.62894$\\
$4$ & $T_4$ & $939039/10240000\approx 0.091703$ & $0.55029$ & $2.62560$\\
$4$ & $\widetilde T_4$ & $\approx 0.091858$ & $0.55052$ & $2.62468$\\
$5$ & $T_5$ & $12654278543/248832000000\approx 0.050855$ & $0.55114$ & $2.62221$\\
$6$ & $T_6$ & $3509759383/124416000000\approx 0.028210$ & $0.55173$ & $2.61984$\\
$1000$ & $T_{1000}$ & $\approx 1.8708\cdot 10^{-256}$ & $0.55497$ & $2.60693$\\
$d\to\infty$ & $T_d$ & & $0.554993$ & $2.60685$\\
\bottomrule
\end{tabular}
\caption{Reducible subsets of tori and the constants they give: $c$ for $(\Zm)^n$
(Corollary~\ref{cor:torus}) and $C=2\sqrt{\log_2(1/c^2)}$ for $r_3(N)$ (Theorem~\ref{thm:behrend}). The values
of $c$ are rounded down and those of $C$ are rounded up.}
\label{tab:constants}
\end{table}

The results from high floors are asymptotic statements: the constant in the exponent $n-O(\log n)$ grows with
the dimension of the reducible set. For a given prime and moderate $n$ the sets in dimensions $2$, $3$ and $4$
are the useful ones, and the volumes of $R_{1/10}$, $T_3$ and $T_4$ can be checked by hand.

In Section~\ref{sec:upper} we prove upper bounds, which limit what such constructions can give. Every
measurable reducible subset of $[0,1)^d$ has volume at most $\frac12$, at most $\frac25$ for $d=2$, and, by the
bound of Ellenberg and Gijswijt, at most $J_\infty^{\,d}$. So no reducible set gives a constant $c$ above
$J_\infty<0.8415$, and the method of~\cite{EHPS} cannot give a constant $C$ below $1.41$. For towers in which
every floor has the floor two dimensions lower as its free set, as in ours, the constants are asymptotically
at most $(1+\sqrt3)/4\approx0.683$, whatever the parameters. These bounds say nothing about $r_3(N)$ itself.

\subsection*{Organisation}
Section~\ref{sec:prelim} contains the facts that we need from~\cite{EPS}. Section~\ref{sec:step} contains
blocks and the induction step. In Section~\ref{sec:tower} we build the tower, compute volumes and prove
Theorem~\ref{thm:main}. Section~\ref{sec:dim3} contains the sets of Theorem~\ref{thm:dim3}. In
Section~\ref{sec:behrend} we prove the induction step for potentials and Theorem~\ref{thm:behrend}.
Section~\ref{sec:upper} contains upper bounds, and Section~\ref{sec:concluding} describes how the constructions
were found and lists open problems. In Appendix~\ref{app:rule} we show that a
seemingly more general removal rule reduces the same sets, and we give a short proof of the lifting principle
for all moduli. The programmes for the computations on which results depend are supplied with the paper.

\section{Reducibility and the lifting principle}\label{sec:prelim}

\subsection*{Notation}
We identify the torus $[0,1)^d$ with the group $\mathbb{R}^d/\Z^d$. Congruences $x\equiv y\pmod 1$
between points of $\mathbb{R}^d$ are meant coordinatewise, and $\vol$ denotes Lebesgue measure. For a real
number $x$ we write $\hf{x}$ for the unique number in $[0,\frac12)$ with $x-\hf{x}\in\frac12\Z$. Then
$2x\equiv 2\hf{x}\pmod 1$. For a point $x\in[0,1)$ the numbers $x$ and $\hf{x}$ agree if and only if
$x<\frac12$; we say that $x$ lies in the \emph{lower half} or the \emph{upper half} of $[0,1)$ according to
whether $x<\frac12$ or not.

\subsection*{Reducible sets}
We use the notion of reducibility from~\cite[Section~2]{EPS}, formulated for arbitrary abelian groups.

\begin{definition}\label{def:reducible}
Let $G$ be an abelian group and $S\subseteq G$. A point $y\in S$ is a \emph{mid-point} of $S$ if there are
distinct $x,z\in S$ with $x+z=2y$, and a \emph{non-mid-point} of $S$ otherwise. The set $S$ is
\emph{non-mid-point reducible}, or \emph{reducible} for short, if every non-empty finite subset of $S$ has a
non-mid-point.
\end{definition}

If $x+z=2y$ with $x\ne z$, then $y\ne x$ and $y\ne z$. A finite set $S$ is reducible if and only if it can be
reduced to the empty set by repeatedly removing non-mid-points~\cite[Section~2]{EPS}, that is, if its points
can be listed as $s_1,\dots,s_N$ so that each $s_i$ is a non-mid-point of $\{s_i,\dots,s_N\}$; we call such a list
a \emph{removal order}. Subsets of reducible sets are reducible. If $R\subseteq[0,1)^d$ is
reducible and $t\in[0,1)$, then the slice $\{x'\in[0,1)^{d-1}:(x',t)\in R\}$ is reducible, because three
points with the same last coordinate satisfy $2y\equiv x+z\pmod 1$ if and only if their first $d-1$
coordinates do.

The following is the main tool of~\cite{EPS}. It turns a reducible set in a fixed dimension into
progression-free sets in all dimensions.

\begin{proposition}[{\cite[Proposition~2.4]{EPS}}]\label{prop:eps}
Let $m\ge 2$ and $d\ge 1$ be fixed, and let $S\subseteq(\Zm)^d$ be non-mid-point reducible. Then for
$n\to\infty$ there is a progression-free subset of $(\Zm)^n$ of size at least $|S|^{\,n/d-O(\log n)}$.
\end{proposition}

The progression-free set consists of the $k$-tuples of elements of $S$, read as elements of
$(\Zm)^{dk}$, in which every element of $S$ occurs equally often. We give a short proof of the
progression-freeness of this set, for odd and for even $m$, in Remark~\ref{rem:eps24}. The exponent
$O(\log n)$ hides a constant that depends on $|S|$, hence on $m$ and $d$; for a set $S$ in dimension~$200$ or
$1000$ it is very large.

\subsection*{\texorpdfstring{From the torus to $(\Zm)^d$}{From the torus to (Z/mZ)\^{}d}}
As in~\cite[Section~3]{EPS}, we construct reducible sets in $(\Zm)^d$ for all $m$ at once, as
discretisations of reducible subsets of the torus. For $m\ge 2$ and $\alpha\in\mathbb{R}^d$ let
$\phi_\alpha\colon(\Zm)^d\to\mathbb{R}^d/\Z^d$ be the map $q\mapsto\alpha+q/m$, and for
$R\subseteq[0,1)^d$ let $S_\alpha(R)=\phi_\alpha^{-1}(R)$.

\begin{lemma}[{cf.~\cite[Section~3]{EPS}}]\label{lem:transfer}
If $R\subseteq[0,1)^d$ is reducible, then $S_\alpha(R)$ is reducible for every $m\ge 2$ and every
$\alpha\in\mathbb{R}^d$.
\end{lemma}

\begin{proof}
The map $\phi_\alpha$ is well defined and injective, and
$2\phi_\alpha(y)-\phi_\alpha(x)-\phi_\alpha(z)\equiv(2y-x-z)/m\pmod 1$. This is $\equiv 0$ if and only if
$2y=x+z$ in $(\Zm)^d$. So $\phi_\alpha$ maps every finite subset of $S_\alpha(R)$ bijectively onto a finite
subset of $R$, and it maps non-mid-points to non-mid-points.
\end{proof}

\begin{lemma}[{cf.~\cite[Section~3]{EPS}}]\label{lem:average}
If $R\subseteq[0,1)^d$ is measurable, then there is $\alpha\in[0,1)^d$ with
$|S_\alpha(R)|\ge\vol(R)\,m^d$.
\end{lemma}

\begin{proof}
If $\alpha$ is uniformly distributed in $[0,1)^d$, then so is $\phi_\alpha(q)$ for every $q$. Hence the
expected size of $S_\alpha(R)$ is $\vol(R)\,m^d$, and some $\alpha$ attains at least the expectation.
\end{proof}

\begin{corollary}\label{cor:torus}
Let $R\subseteq[0,1)^d$ be a measurable reducible set of volume $v$, and let $m\ge 2$ be such that
$v^{1/d}m>1$. Then for $n\to\infty$ there is a progression-free subset of $(\Zm)^n$ of size at least
$(v^{1/d}m)^{\,n-O(\log n)}$, where the implied constant depends on $d$, $m$ and $v$.
\end{corollary}

\begin{proof}
By Lemmas~\ref{lem:transfer} and~\ref{lem:average} there is a reducible set $S\subseteq(\Zm)^d$ with
$vm^d\le|S|\le m^d$. Proposition~\ref{prop:eps} gives a progression-free set of size at least
$|S|^{n/d}\,|S|^{-O(\log n)}$, which is at least $(v^{1/d}m)^n\,m^{-O(d\log n)}$. Since $v^{1/d}m>1$, the
last factor is $(v^{1/d}m)^{-O(\log n)}$.
\end{proof}

All the seeds used for our lower bounds have $v^{1/d}\ge\frac12$, so the condition $v^{1/d}m>1$ holds for all $m\ge 2$
except for $[0,\frac12)$ with $m=2$. In view of Corollary~\ref{cor:torus}, all constructions below take
place in the torus.

\subsection*{An elementary fact}
All our proofs of reducibility rest on the following observation about mid-points of numbers whose
reductions modulo $\frac12$ lie in a short interval.

\begin{lemma}\label{lem:interval}
Let $I=[a,a+\ell)$ be an interval with $\ell\le\frac14$ contained in $[0,\frac12)$, and let
$x,y,z\in[0,1)$ satisfy $\hf{x},\hf{y},\hf{z}\in I$ and $2y\equiv x+z\pmod 1$. Then $x$ and $z$ lie in the
same half of $[0,1)$, and $\hf{y}=\frac12(\hf{x}+\hf{z})$. If instead $x,y,z\in[0,\frac12)$ and
$2y\equiv x+z\pmod 1$, then $y=\frac12(x+z)$.
\end{lemma}

\begin{proof}
Write $x=\hf{x}+\frac i2$ and $z=\hf{z}+\frac k2$ with $i,k\in\{0,1\}$. The numbers $2\hf{y}-2a$ and
$\hf{x}+\hf{z}-2a$ lie in $[0,2\ell)\subseteq[0,\frac12)$, so their difference
$D=\hf{x}+\hf{z}-2\hf{y}$ lies in $(-\frac12,\frac12)$. Since $2y\equiv2\hf{y}$, the congruence says that
$D+\frac{i+k}2$ is an integer. If $i+k=1$, then $D\equiv\frac12\pmod 1$, which is impossible. Hence $i=k$,
so $D$ is an integer, and $D=0$. For the second statement, $2y$ and $x+z$ both lie in $[0,1)$, and two
congruent numbers in $[0,1)$ are equal.
\end{proof}

The next lemma is the two-dimensional version that we need for the low zones of our sets, in which rows
are arranged in a window that slides with the level.

\begin{lemma}\label{lem:chart}
Let $0\le\beta\le\frac14$ and
\[
  \Omega=\Omega_\beta=\bigl\{(\sigma,\tau):\ 0\le\tau<\tfrac14,\ \ \max(0,\tau-\tfrac14+\beta)\le\sigma<\tau+\beta\bigr\}.
\]
Let $(s_x,t_x)$, $(s_y,t_y)$, $(s_z,t_z)$ be points of $[0,1)^2$ with $s_x+s_z\equiv2s_y$ and
$t_x+t_z\equiv2t_y\pmod1$, and suppose that the points $(\sigma,\tau)=(\hf{s},\hf{t})$ of all three lie in
$\Omega$. Then $\sigma_x+\sigma_z=2\sigma_y$ and $\tau_x+\tau_z=2\tau_y$. Moreover $s_x$ and $s_z$ lie in the same
half of $[0,1)$, and so do $t_x$ and $t_z$.
\end{lemma}

\begin{proof}
The numbers $\sigma_x+\sigma_z-2\sigma_y$ and $\tau_x+\tau_z-2\tau_y$ are multiples of $\frac12$. All values of
$\tau$ lie in $[0,\frac14)$, so $\tau_x+\tau_z=2\tau_y$. The values of $\sigma$ lie in $[0,\frac12)$, so
$\sigma_x+\sigma_z-2\sigma_y$ is $0$ or $\pm\frac12$. Let $\phi(\tau)=\max(0,\tau-\frac14+\beta)$, a convex
function with $\phi(\tau)+\frac14\ge\tau+\beta$. If $\sigma_x+\sigma_z=2\sigma_y+\frac12$, then
$2\sigma_y+\frac12<(\tau_x+\beta)+(\tau_z+\beta)=2\tau_y+2\beta$, so $\sigma_y<\tau_y+\beta-\frac14\le\phi(\tau_y)$.
If $\sigma_x+\sigma_z=2\sigma_y-\frac12$, then $2\sigma_y-\frac12\ge\phi(\tau_x)+\phi(\tau_z)\ge2\phi(\tau_y)$, so
$\sigma_y\ge\phi(\tau_y)+\frac14\ge\tau_y+\beta$. Both contradict $(\sigma_y,\tau_y)\in\Omega$, so
$\sigma_x+\sigma_z=2\sigma_y$. Finally $s_x-\sigma_x$ and $s_z-\sigma_z$ lie in $\{0,\frac12\}$, and their sum is
congruent to $2s_y-2\sigma_y\equiv0$ modulo~$1$, so they are equal; likewise for $t$.
\end{proof}

\section{Blocks and the induction step}\label{sec:step}

We write a point of $[0,1)^d$ as $(x',t)$ with $x'\in[0,1)^{d-1}$, and call $t$ its \emph{level}. For a
set $B\subseteq[0,1)^d$ and $t\in[0,1)$ the \emph{slice} of $B$ at level $t$ is
$B_t=\{x'\in[0,1)^{d-1}:(x',t)\in B\}$. For $d=1$ the set $[0,1)^0$ is a single point, and every slice is
either empty or this point.

\begin{definition}\label{def:block}
Let $d\ge 1$. A set $B\subseteq[0,1)^d$ is a \emph{block} if it has the following three properties.
\begin{itemize}
\item[(L)] Every point of $B$ has level $t<\frac34$.
\item[(R)] $B$ is non-mid-point reducible.
\item[(B)] If $x,y,z\in B$ have levels with $t_x+t_z=2t_y+1$, then $x'+z'\not\equiv 2y'\pmod 1$.
\end{itemize}
\end{definition}

\begin{remark}\label{rem:levels}
In property~(B) the levels satisfy $t_x,t_z>\frac14$ and $t_y<\frac14$. Indeed, $t_x+t_z\ge 1$ and
$t_x,t_z<\frac34$ force $t_x,t_z>\frac14$, and then $t_y=\frac12(t_x+t_z-1)<\frac14$. So property~(B)
concerns triples that wrap around the circle once in the level coordinate, with the two outer points
high and the middle point low. Triples with $t_x+t_z=2t_y-1$ need not be excluded: they never occur in a
removal order that starts at the lowest level. For $d=1$, property~(B) says that $B$ contains no three
levels with $t_x+t_z=2t_y+1$.
\end{remark}

In the terminology of~\cite{EPS}, every reducible set is a building block for progression-free sets in high
dimensions; a block is a building block with the two additional properties~(L) and~(B) of its last coordinate.
Given these two properties, reducibility only has to be checked inside the slices.

\begin{lemma}\label{lem:blockchar}
Let $B\subseteq[0,1)^d$ have properties~(L) and~(B). Then $B$ is a block if and only if every slice $B_t$ is
non-mid-point reducible.
\end{lemma}

\begin{proof}
Slices of reducible sets are reducible (Section~\ref{sec:prelim}). Conversely, let every slice be reducible,
let $S$ be a non-empty finite subset of $B$, let $t_0$ be the lowest level of a point of $S$, and let $S_0$ be
the set of points of $S$ at level $t_0$. Since the slice at level $t_0$ is reducible, $S_0$ has a
non-mid-point $y$. Suppose that $2y\equiv x+z\pmod 1$ with distinct $x,z\in S$. Then $t_x+t_z-2t_0$ is an
integer in $[0,\frac32)$, so it is $0$ or $1$. If it is $0$, then $t_x=t_z=t_0$, so $x$ and $z$ lie in
$S_0$, against the choice of $y$. If it is $1$, then $t_x+t_z=2t_y+1$, and property~(B) shows that
$2y\not\equiv x+z$. Hence $y$ is a non-mid-point of $S$.
\end{proof}

The induction step takes a block $B$ in dimension $d$ and a reducible set $F$ in dimension $d-1$, and
produces a set in dimension $d+1$. We write its points as $(x',t,w)$ with $x'\in[0,1)^{d-1}$ and
$t,w\in[0,1)$. The new coordinate $w$ is the level of the new set; the old level $t$ becomes an ordinary
coordinate.

\begin{definition}\label{def:step}
Let $B\subseteq[0,1)^d$ be a block, let $F\subseteq[0,1)^{d-1}$ be reducible, and let
$0\le\beta\le\frac14$. For $0\le\tau<\frac12$ let the \emph{row set} $P_\tau$ be $B_\tau$ if
$\tau<\frac14$ and $F$ if $\tau\ge\frac14$. The set $B^+=B^+(F,\beta)\subseteq[0,1)^{d+1}$ is the union of
the \emph{high zone}
\[
  H=\bigl\{(x',t,w):\ \tfrac14\le w<\tfrac34,\ \ (x',t)\in B,\ \ t\ge w-\tfrac12+\beta\bigr\}
\]
and the \emph{low zone}
\[
  L=\bigl\{(x',t,w):\ 0\le w<\tfrac14,\ \ \max(0,\,w-\tfrac14+\beta)\le\hf{t}<w+\beta,\ \
  x'\in P_{\hf{t}}\bigr\}.
\]
\end{definition}

\begin{center}
\small
\renewcommand{\arraystretch}{1.3}
\begin{tabular}{@{}lll@{}}
\toprule
Part of $B^+$ & Levels $w$ & Condition on $(x',t)$\\
\midrule
High zone & $\frac14\le w<\frac34$ & $(x',t)\in B$ and $t\ge w-\frac12+\beta$\\
Low zone, inherited rows & $0\le w<\frac14$ & $\max(0,w-\frac14+\beta)\le\hf{t}<\min(\frac14,w+\beta)$ and $x'\in B_{\hf{t}}$\\
Low zone, free rows & $\frac14-\beta<w<\frac14$ & $\frac14\le\hf{t}<w+\beta$ and $x'\in F$\\
\bottomrule
\end{tabular}
\end{center}

\begin{figure}[tbp]
\centering
\includegraphics[width=0.78\textwidth]{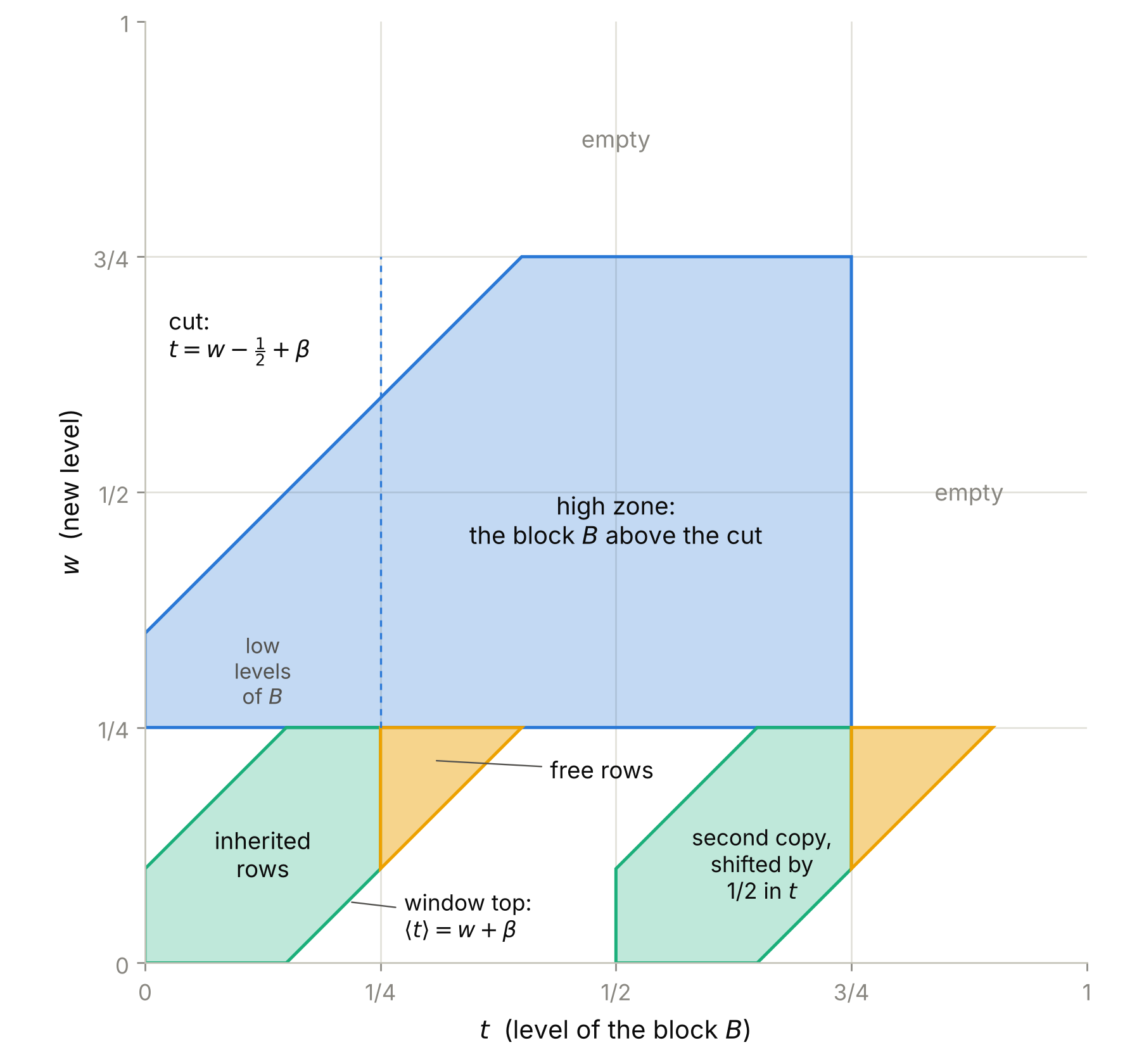}
\caption{The induction step in the plane of the old level $t$ and the new level $w$, drawn for
$\beta=3/20$. The high zone is the part of $B$ above the cut $t=w-\frac12+\beta$. The low zone is the
window $\max(0,w-\frac14+\beta)\le\hf{t}<w+\beta$, taken twice. A row with $\hf{t}<\frac14$ carries the
slice of $B$ at level $\hf{t}$, and a row with $\hf{t}\ge\frac14$ carries the free set $F$.}
\label{fig:step}
\end{figure}

In the low zone only $\hf{t}$ is restricted, so $t$ may be $\hf{t}$ or $\hf{t}+\frac12$, and every row
occurs twice. For each $w<\frac14$ the values of $\hf{t}$ form the \emph{window}
$I_w=[\max(0,w-\frac14+\beta),\,w+\beta)$, a half-open interval of length at most $\frac14$ inside
$[0,\frac12)$, which slides upwards as $w$ grows. Figure~\ref{fig:step} shows the two zones in the plane
of the levels $t$ and $w$.

The proof of the following theorem rests on two ideas. First, the exclusion of wrapped triples is inherited:
a wrapped triple of $B^+$ has its middle point in an inherited row and projects to a wrapped triple of $B$,
which property~(B) of $B$ excludes. Second, rows in a window of length at most $\frac14$ can be used twice, at
$t$ and at $t+\frac12$, by Lemma~\ref{lem:interval}; this is what can make the low slices larger than those of a
product.

\begin{theorem}\label{thm:step}
If $B$ is a block, then $B^+$ is a block, with $w$ as its level. This holds for every $0\le\beta\le\frac14$ and
every reducible set $F$.
\end{theorem}

\begin{proof}
Property~(L) holds because $w<\frac34$ on $B^+$.

\emph{Property~(B).} Let $x,y,z$ be points of $B^+$ with $w_x+w_z=2w_y+1$. By Remark~\ref{rem:levels},
$w_x,w_z>\frac14$ and $w_y<\frac14$, so $x$ and $z$ lie in the high zone and $y$ lies in the low zone. Put
$\tau=\hf{t_y}$, so that $2t_y\equiv2\tau\pmod 1$. The window gives $\tau<w_y+\beta$, the cut of the high
zone gives $t_x\ge w_x-\frac12+\beta$ and $t_z\ge w_z-\frac12+\beta$, and property~(L) of $B$ gives
$t_x,t_z<\frac34$. Hence
\[
  2\tau\;<\;2w_y+2\beta\;=\;(w_x-\tfrac12+\beta)+(w_z-\tfrac12+\beta)\;\le\;t_x+t_z\;<\;\tfrac32 .
\]
Suppose that $t_x+t_z\equiv 2t_y\pmod 1$. Then $t_x+t_z-2\tau$ is an integer in $(0,\frac32)$, so
$t_x+t_z=2\tau+1$. Since $t_x+t_z<\frac32$, this gives $\tau<\frac14$. So the row of $y$ is inherited:
$y'\in B_\tau$, that is, $(y',\tau)\in B$. The points $(x',t_x)$ and $(z',t_z)$ lie in $B$ as well, and the
levels of the three points $(x',t_x)$, $(y',\tau)$, $(z',t_z)$ of $B$ satisfy $t_x+t_z=2\tau+1$. By
property~(B) of $B$, $x'+z'\not\equiv2y'\pmod 1$. So the congruence $2y\equiv x+z$ fails in one of the
first $d$ coordinates, which proves property~(B) for $B^+$.

\emph{Property~(R).} By Lemma~\ref{lem:blockchar} it suffices to show that every slice of $B^+$ is reducible. For $w\ge\frac14$ the slice is a
subset of $B$. For $w<\frac14$ the slice consists of the points $(x',t)$ with $\hf{t}\in I_w$ and
$x'\in P_{\hf{t}}$, and every row set $P_\tau$ is reducible, as a slice of $B$ or as $F$. Let $G$ be a
non-empty finite subset of this slice. Let $\tau_0$ be the largest value of $\hf{t}$ among the points of
$G$, and let $P$ be the set of all $x'$ such that $G$ contains a point $(x',t)$ with $\hf{t}=\tau_0$. Then
$P$ is a non-empty finite subset of $P_{\tau_0}$, so it has a non-mid-point $p$. Choose a point $y\in G$
with $\hf{t_y}=\tau_0$ and $y'=p$. Suppose that $2y\equiv x+z\pmod 1$ with distinct $x,z\in G$. By
Lemma~\ref{lem:interval}, applied to the coordinate $t$ and the interval $I_w$, the numbers $t_x$ and
$t_z$ lie in the same half of $[0,1)$ and $\tau_0=\frac12(\hf{t_x}+\hf{t_z})$. Since $\tau_0$ is the largest
value, $\hf{t_x}=\hf{t_z}=\tau_0$, and hence $t_x=t_z$ and $x'\ne z'$. Both $x'$ and $z'$ lie in $P$, and
$x'+z'\equiv2p\pmod 1$, which contradicts the choice of $p$. So $y$ is a non-mid-point of $G$.
\end{proof}

The proof gives an explicit removal order for $B^+$: the levels are removed from the lowest upwards, and
inside a level of the low zone the points are removed by decreasing $\hf{t}$, in each row in a removal
order of the row set.

\begin{remark}\label{rem:step}
The proof uses only a few features of the definition.
\begin{enumerate}[label=(\alph*)]
\item The free rows never take part in the wrapped triples that property~(B) excludes: in the proof of
property~(B) the middle point is forced into a row with $\hf{t}<\frac14$. The free rows therefore only need to be reducible, and $F$ can be
any reducible set.
\item The upper end $w+\beta$ of the window and the cut $w-\frac12+\beta$ of the high zone can be replaced
by functions $U(w)$ and $C(w)$, as long as $C(w_x)+C(w_z)\ge 2U(w_y)$ whenever $w_x+w_z=2w_y+1$, and as long as
the other features are kept: the windows are half-open intervals of length at most $\frac14$ inside
$[0,\frac12)$, the rows with $\hf{t}<\frac14$ are inherited from $B$, and the other rows are reducible. The
lower end of the window only serves to keep the window shorter than $\frac14$, which the slices need. If the
cut $C$ is convex, then $U(w)=C(w+\frac12)$ satisfies the condition, by Jensen's inequality, so a convex cut
determines an admissible window. In numerical experiments in dimension~$3$, curved cuts and windows changed
the optimal volume by less than $10^{-4}$.
\item The step can gain over the product of $B$ with an interval of length $\frac12$ because it uses all
levels $w\in[0,\frac34)$, with different slices at different levels. At the levels $w<\frac14$ every row is
used twice, at $t=\hf{t}$ and at $t=\hf{t}+\frac12$. It does not gain for every block: for $B=[0,\varepsilon)$
with small $\varepsilon$, $F=\emptyset$ and $\beta=\frac14$, formula~\eqref{eq:volplus} below gives
$\vol B^+=\frac32\varepsilon^2<\frac12\varepsilon$. For the blocks of the tower it does gain.
\end{enumerate}
\end{remark}

\section{The tower: volumes and growth rate}\label{sec:tower}

Throughout this section $\frac14\le\gamma\le\frac12$. We use the two profile functions
\[
  \rho(r)=\min\bigl(\tfrac14,\;r+\tfrac12-\gamma\bigr), \qquad \chi(s)=\max\bigl(0,\;s-\gamma\bigr).
\]

\subsection*{The first floor: a block in dimension 2}

\begin{definition}\label{def:E}
Let $E(\gamma)\subseteq[0,1)^2$ be the set of points $(u,s)$ with
\[
  0\le s<\tfrac14\ \text{ and }\ \hf{u}<\rho(s), \qquad\text{or}\qquad
  \tfrac14\le s<\tfrac34\ \text{ and }\ \chi(s)\le u<\tfrac12 .
\]
We call the first part the \emph{doubled part}, since it consists of two copies that differ by
$\frac12$ in $u$, and the second part the \emph{single part}. The level of $E(\gamma)$ is $s$.
\end{definition}

\begin{lemma}\label{lem:E}
$E(\gamma)$ is a block.
\end{lemma}

\begin{proof}
Property~(L) is clear.

\emph{Property~(B).} Let $x,y,z\in E(\gamma)$ with $s_x+s_z=2s_y+1$. By Remark~\ref{rem:levels},
$s_x,s_z>\frac14$ and $s_y<\frac14$, so $x$ and $z$ lie in the single part and $y$ lies in the doubled
part. Hence
\[
  u_x+u_z\;\ge\;\chi(s_x)+\chi(s_z)\;\ge\;s_x+s_z-2\gamma\;=\;2s_y+1-2\gamma\;>\;2\hf{u_y},
\]
because $\hf{u_y}<\rho(s_y)\le s_y+\frac12-\gamma$. Also $u_x+u_z<1$. The two numbers $2\hf{u_y}$ and
$u_x+u_z$ lie in $[0,1)$ and are different, so they are not congruent modulo~$1$. Since
$2u_y\equiv2\hf{u_y}$, we get $u_x+u_z\not\equiv2u_y\pmod 1$.

\emph{Property~(R).} By Lemma~\ref{lem:blockchar} it suffices to show that every slice is reducible. In a
non-empty finite subset of a slice, a point $y$ with the largest value of $\hf{u}$ is a non-mid-point: if
$2y\equiv x+z$ with distinct $x,z$ in the subset, then by Lemma~\ref{lem:interval}, with $I=[0,\rho(s))$ for a
level $s<\frac14$ and by its second part for $s\ge\frac14$, the number $\hf{u_y}$ is the average of the two
different numbers $\hf{u_x}$ and $\hf{u_z}$, so it is smaller than one of them.
\end{proof}

\begin{remark}\label{rem:EPS}
The area of $E(\gamma)$ is computed below; it is $7/24$ for $\gamma=5/12$, which is
the largest value, and $233/800$ for $\gamma=2/5$. The set $E(5/12)$ is closely related to the
two-dimensional set of~\cite[Section~3]{EPS}, which also has area $7/24$ and which is supported on the points
$(a,b)$ with $\frac7{12}\le a+b\le\frac43$. In the coordinates $u=a$ and $s=a+b-\frac7{12}$ modulo~$1$, an
affine change of coordinates of the torus with unimodular linear part, the two sets differ, up to their
boundaries (the set of~\cite{EPS} is closed, ours is half-open), only in the rows $\frac16\le s<\frac14$, where
the fibre $[0,s+\frac1{12}]\cup[\frac12,\frac{11}{12}-s)$ is replaced by the two arcs $\hf{u}<\frac14$ of the
same total length.
\end{remark}

\subsection*{The second floor: a family in dimension 3}

Applying the induction step to $E(\gamma)$ gives the following family of sets in dimension~$3$, with
coordinates $(u,s,t)$ and level $t$.

\begin{definition}\label{def:T}
For $0\le\beta\le\frac14$ let $T(\beta,\gamma)\subseteq[0,1)^3$ be the union of the \emph{low zone}
\[
  \bigl\{(u,s,t):\ 0\le t<\tfrac14,\ \ \max(0,\,t-\tfrac14+\beta)\le\hf{s}<t+\beta,\ \
  \hf{u}<\rho(\hf{s})\bigr\}
\]
and the \emph{high zone}
\[
  \bigl\{(u,s,t):\ \tfrac14\le t<\tfrac34,\ \ (u,s)\in E(\gamma),\ \ s\ge t-\tfrac12+\beta\bigr\}.
\]
\end{definition}

\begin{center}
\small
\renewcommand{\arraystretch}{1.3}
\begin{tabular}{@{}llllc@{}}
\toprule
Part of $T(\beta,\gamma)$ & Levels $t$ & Range of $s$ & Range of $u$ & Copies\\
\midrule
Low zone & $0\le t<\frac14$ & $\max(0,t-\frac14+\beta)\le\hf{s}<t+\beta$ & $\hf{u}<\rho(\hf{s})$ & $4$\\
Doubled part & $\frac14\le t<\frac34$ & $\max(0,t-\frac12+\beta)\le s<\frac14$ & $\hf{u}<\rho(s)$ & $2$\\
Single part & $\frac14\le t<\frac34$ & $\max(\frac14,t-\frac12+\beta)\le s<\frac34$ & $\chi(s)\le u<\frac12$ & $1$\\
\bottomrule
\end{tabular}
\end{center}

The column ``Copies'' counts the translates by $\frac12$ in $u$ and in $s$: the low zone consists of four
copies of one region, the doubled part of two. At the levels $\frac14\le t\le\frac12-\beta$ the slice is
all of $E(\gamma)$; above them it is cut back from below in $s$ by the line $s=t-\frac12+\beta$. At the
levels $t<\frac14$ the slice lies under the graph of $\rho$ over a window of values of $\hf{s}$ of length
at most $\frac14$, which slides upwards with $t$.

\begin{proposition}\label{prop:TEplus}
Let $F_1=\{u\in[0,1):\hf{u}<\frac14\}$. Then $F_1$ is reducible, and $T(\beta,\gamma)=E(\gamma)^+(F_1,\beta)$,
where the coordinates $(u,s,t)$ of $T(\beta,\gamma)$ take the roles of $(x',t,w)$ in
Definition~\ref{def:step}.
\end{proposition}

\begin{proof}
In a non-empty finite subset of $F_1$, a point $y$ with the largest value of $\hf{u}$ is a non-mid-point:
by Lemma~\ref{lem:interval} with $I=[0,\frac14)$, a solution of $2y\equiv x+z$ would have $x$ and $z$ in the
same half and $\hf{y}$ equal to the average of $\hf{x}$ and $\hf{z}$, so $\hf{x}=\hf{z}=\hf{y}$ and $x=z$.
The high zones of the two sets agree by definition. In the low zone, a row with $\hf{s}<\frac14$ is the
slice of $E(\gamma)$ at level $\hf{s}$, which is $\{\hf{u}<\rho(\hf{s})\}$. A row with $\hf{s}\ge\frac14$ is
$F_1$, and $\rho(\hf{s})=\frac14$ for these rows, because $\hf{s}+\frac12-\gamma\ge\frac14$.
\end{proof}

\begin{corollary}\label{cor:T}
For $0\le\beta\le\frac14$ and $\frac14\le\gamma\le\frac12$ the set $T(\beta,\gamma)$ is a block. In
particular it is non-mid-point reducible.
\end{corollary}

\begin{proof}
This follows from Lemma~\ref{lem:E}, Proposition~\ref{prop:TEplus} and Theorem~\ref{thm:step}.
\end{proof}

Unwinding the proofs gives an explicit removal order for $T(\beta,\gamma)$: lowest $t$ first; at a level
$t<\frac14$, by decreasing $\hf{s}$ and then by decreasing $\hf{u}$; at a level $t\ge\frac14$, by
increasing $s$ and then by decreasing $\hf{u}$.

\subsection*{All floors}

\begin{definition}\label{def:tower}
Let $\gamma\in[\frac14,\frac12]$ and $\beta_3,\beta_4,\ldots\in[0,\frac14]$. The \emph{tower} is the
sequence of sets $T_d\subseteq[0,1)^d$ defined by
\[
  T_2=E(\gamma),\qquad T_3=T(\beta_3,\gamma),\qquad T_4=T_3^+\bigl(E(\tfrac5{12}),\beta_4\bigr),\qquad
  T_{d+1}=T_d^+\bigl(T_{d-1},\beta_{d+1}\bigr)\quad(d\ge 4).
\]
Unless stated otherwise we use $\gamma=\frac25$, $\beta_3=\frac18$, $\beta_4=\frac3{20}$ and
$\beta_d=\frac16$ for $d\ge5$.
\end{definition}

So the free rows of $T_3$ are copies of $F_1$, of measure $\frac12$, the free rows of $T_4$ are copies of
the largest set $E(5/12)$ of the first floor, and from then on the free rows of $T_{d+1}$ are copies of
the floor $T_{d-1}$ two dimensions lower.

\begin{theorem}\label{thm:tower}
For every $d\ge2$ the set $T_d$ is a block. In particular it is non-mid-point reducible.
\end{theorem}

\begin{proof}
By induction on $d$. The sets $T_2$ and $T_3$ are blocks by Lemma~\ref{lem:E} and
Corollary~\ref{cor:T}. The free sets $E(5/12)$ and $T_{d-1}$ are reducible by Lemma~\ref{lem:E} and the
induction hypothesis, so Theorem~\ref{thm:step} shows that $T_{d+1}$ is a block.
\end{proof}

In coordinates $(u,s,t,w)$, with $T_3=T(\frac18,\frac25)$ and $\beta=\beta_4$, the fourth floor is the
following set.

\begin{center}
\small
\renewcommand{\arraystretch}{1.3}
\begin{tabular}{@{}ll>{\raggedright\arraybackslash}p{7.6cm}@{}}
\toprule
Part of $T_4$ & Levels $w$ & Condition on $(u,s,t)$\\
\midrule
High zone & $\frac14\le w<\frac34$ & $(u,s,t)\in T_3$ and $t\ge w-\frac12+\beta$\\
Inherited rows & $0\le w<\frac14$ & $\max(0,w-\frac14+\beta)\le\hf{t}<\min(\frac14,w+\beta)$,
$\max(0,\hf{t}-\frac18)\le\hf{s}<\hf{t}+\frac18$, $\hf{u}<\rho(\hf{s})$\\
Free rows & $\frac14-\beta<w<\frac14$ & $\frac14\le\hf{t}<w+\beta$ and $(u,s)\in E(\frac5{12})$\\
\bottomrule
\end{tabular}
\end{center}

\begin{remark}\label{rem:floor1}
The tower can be started one floor lower. The interval $[0,\frac12)\subseteq[0,1)$ is a block in
dimension~$1$: its levels are below $\frac34$, a largest element of a finite subset is a non-mid-point by
the second part of Lemma~\ref{lem:interval}, and property~(B) holds because $t_x+t_z<1$. Its
discretisations are the intervals behind the constant $\frac12$, as in the construction of Salem and
Spencer~\cite{SalemSpencer}. Apply the step to it with $F=[0,1)^0$, a single point, and with
$\beta=\frac12-\gamma$, and write $(u,s)$ for the coordinates $(t,w)$. The high zone is the single part
of $E(\gamma)$, because $w-\frac12+\beta=s-\gamma$. In the low zone the rows with $s\le\gamma-\frac14$ are
$\hf{u}<s+\beta=\rho(s)$, as in $E(\gamma)$. For $\gamma-\frac14<s<\frac14$ the row is the arc
$s-\frac14+\beta\le\hf{u}<s+\beta$ of length $\frac14$ instead of $\hf{u}<\frac14$. So
$[0,\frac12)^+$ is obtained from $E(\gamma)$ by rotating the arcs in these rows. It has the same slice
measures as $E(\gamma)$, so a tower started from it has the same volumes (Proposition~\ref{prop:volstep}).
In this sense the construction of~\cite{EPS} is the first step of the induction. In Section~\ref{sec:dim3} we add
a notch to $[0,\frac12)^+$ and obtain a larger reducible set, which is not a block.
\end{remark}

\subsection*{The volume recursion}

For a Borel set $B\subseteq[0,1)^d$ let $a_B(t)$ denote the $(d-1)$-dimensional measure of the
slice $B_t$; we call $a_B$ the \emph{slice measure} of $B$. All sets that we construct are Borel sets, so all
their slices are measurable. (For Lebesgue measurable sets the statements below hold for almost every level.) The volume of $B^+$ depends on $B$ only through
$a_B$ and on $F$ only through its volume.

\begin{proposition}\label{prop:volstep}
Let $B$, $F$ and $\beta$ be as in Definition~\ref{def:step}, with $B$ and $F$ Borel sets, and write
$b=a_B$ and $f=\vol F$. Let $r(\tau)=b(\tau)$ for $0\le\tau<\frac14$ and $r(\tau)=f$ for
$\frac14\le\tau<\frac12$. Then the slice measure of $B^+$ is
\begin{align}
  a_{B^+}(w)&=\int_{\max(0,\,w-\frac12+\beta)}^{3/4} b(t)\,dt &&\text{for }\tfrac14\le w<\tfrac34,
  \label{eq:ahigh}\\
  a_{B^+}(w)&=2\int_{\max(0,\,w-\frac14+\beta)}^{w+\beta} r(\tau)\,d\tau &&\text{for }0\le w<\tfrac14,
  \label{eq:alow}
\end{align}
and
\begin{equation}\label{eq:volplus}
  \vol B^+\;=\;\int_0^{3/4} b(t)\,\omega_h(t)\,dt\;+\;2\int_0^{1/4} b(t)\,\omega_\ell(t)\,dt\;+\;\beta^2 f,
\end{equation}
where $\omega_h(t)=\min(\frac12,\,t+\frac14-\beta)$ and $\omega_\ell(t)=\frac14-|t-\beta|$, that is,
$\omega_\ell(t)=t+\frac14-\beta$ for $0\le t<\beta$ and $\omega_\ell(t)=\frac14+\beta-t$ for $\beta\le t<\frac14$.
\end{proposition}

\begin{proof}
At a level $w\ge\frac14$ the slice of $B^+$ is $\{(x',t)\in B:t\ge w-\frac12+\beta\}$, and $b(t)=0$ for
$t\ge\frac34$; this gives~\eqref{eq:ahigh}. At a level $w<\frac14$ the slice consists, for every
$\tau\in I_w$, of the two rows $t=\tau$ and $t=\tau+\frac12$, each of measure $r(\tau)$; this
gives~\eqref{eq:alow}. To obtain~\eqref{eq:volplus} we integrate over $w$ and exchange the order of
integration. A point $(x',t)\in B$ contributes to the high zone for all
$w\in[\frac14,\frac34)$ with $w\le t+\frac12-\beta$, a set of length $\omega_h(t)$. A value
$\tau\in[0,\frac12)$ lies in the window $I_w$ if and only if $\tau-\beta<w\le\tau+\frac14-\beta$; for
$\tau<\frac14$ the set of such $w\in[0,\frac14)$ has length $\omega_\ell(\tau)$, and for $\tau\ge\frac14$ it
has length $\max(0,\frac14+\beta-\tau)$. Since $\int_{1/4}^{1/4+\beta}(\frac14+\beta-\tau)\,d\tau=\frac12\beta^2$,
the free rows contribute $2f\cdot\frac12\beta^2=\beta^2 f$.
\end{proof}

So the slice measures of the floors of the tower satisfy a linear recursion: $a_{T_{d+1}}$ is determined
by $a_{T_d}$ and by the number $\vol T_{d-1}$. If the parameters $\gamma$ and $\beta_d$ are rational, as they
are for the default parameters of Definition~\ref{def:tower}, then all slice measures are piecewise polynomial functions with
rational breakpoints and rational coefficients, so all volumes are rational numbers.

\subsection*{The first two floors}

By Definition~\ref{def:E}, the slice measure of $E(\gamma)$ is $2\rho(s)$ for $s<\frac14$ and
$\frac12-\chi(s)$ for $\frac14\le s<\frac34$, that is,
\[
  a_E(s)=\begin{cases}
  1-2\gamma+2s & \text{for } 0\le s<\gamma-\tfrac14,\\
  \tfrac12 & \text{for } \gamma-\tfrac14\le s<\gamma,\\
  \tfrac12+\gamma-s & \text{for } \gamma\le s<\tfrac34.
  \end{cases}
\]
Integrating gives
\[
  \vol E(\gamma)\;=\;\tfrac18+\bigl(\tfrac34-\gamma\bigr)\bigl(\tfrac32\gamma-\tfrac18\bigr)
  \;=\;\tfrac7{24}-\tfrac32\bigl(\gamma-\tfrac5{12}\bigr)^2 .
\]
This is largest for $\gamma=\frac5{12}$, where it equals $\frac7{24}$; for $\gamma=\frac25$ it is
$\frac{233}{800}$.

\begin{proposition}\label{prop:volT}
The volume of $T(\beta,\gamma)$ is given by~\eqref{eq:volplus} with $b=a_E$ and $f=\frac12$. Explicitly,
$\vol T(\beta,\gamma)=V_1(\beta,\gamma)$ for $0\le\beta\le\gamma-\frac14$ and
$\vol T(\beta,\gamma)=V_2(\beta,\gamma)$ for $\gamma-\frac14\le\beta\le\frac14$, where
\begin{align*}
  V_1&=-\tfrac43\beta^3+\bigl(4\gamma-\tfrac74\bigr)\beta^2+\bigl(\tfrac1{16}+\tfrac12\gamma-\gamma^2\bigr)\beta
       +\tfrac13\gamma^3-\tfrac54\gamma^2+\tfrac{13}{16}\gamma+\tfrac1{96},\\
  V_2&=\tfrac16\beta^3-\bigl(\tfrac12\gamma+\tfrac58\bigr)\beta^2+\bigl(\tfrac{11}{32}-\tfrac74\gamma+\tfrac72\gamma^2\bigr)\beta
       -\tfrac76\gamma^3-\tfrac18\gamma^2+\tfrac{17}{32}\gamma+\tfrac{13}{384}.
\end{align*}
In particular $\vol T(\tfrac18,\tfrac25)=V_1(\tfrac18,\tfrac25)=\tfrac{15803}{96000}\approx0.164615$.
\end{proposition}

\begin{proof}
The first statement is Propositions~\ref{prop:TEplus} and~\ref{prop:volstep}. The breakpoints of $a_E$ are
$\gamma-\frac14$ and $\gamma$, those of the weights are $\beta$ and $\frac14+\beta$, and their order is the same
throughout each of the two ranges of $\beta$. So each term of~\eqref{eq:volplus} is a sum of integrals of
products of two linear functions, and these give $V_1$ and $V_2$.
\end{proof}

The polynomial $V_1$ has a stationary point at $\beta\approx0.12651$, $\gamma\approx0.39927$, with value
$\approx0.1646174$, so the choice $(\frac18,\frac25)$ is within $3\cdot10^{-6}$ of it. We have not proved that
this is the global maximum of the family.

\subsection*{The fourth floor}

From now on $T=T_3=T(\frac18,\frac25)$.

\begin{proposition}\label{prop:volT4}
(a) The slice measure $a_T$ of $T$ is the piecewise polynomial function given in Table~\ref{tab:aT}.

(b) For $\beta_4=\frac3{20}$,
\[
  \vol T_4\;=\;\underbrace{\tfrac{412081}{6144000}}_{\text{high zone}}
  \;+\;\underbrace{\tfrac{69389}{3840000}}_{\text{inherited rows}}
  \;+\;\underbrace{\tfrac{21}{3200}}_{\text{free rows}}
  \;=\;\tfrac{939039}{10240000}\;\approx\;0.0917030 .
\]
\end{proposition}

\begin{table}[ht]
\centering
\small
\renewcommand{\arraystretch}{1.3}
\begin{tabular}{@{}llll@{}}
\toprule
Levels $t$ & $a_T(t)$ & Levels $t$ & $a_T(t)$\\
\midrule
$[0,\frac1{40})$ & $\frac{13}{160}+\frac9{10}t+2t^2$ & $[\frac14,\frac38)$ & $\frac{233}{800}$\\
$[\frac1{40},\frac18)$ & $\frac2{25}+t$ & $[\frac38,\frac{21}{40})$ & $\frac{361}{1600}+\frac{11}{20}t-t^2$\\
$[\frac18,\frac14)$ & $\frac{79}{800}+\frac{11}{10}t-2t^2$ & $[\frac{21}{40},\frac34)$ & $\frac{401}{800}-\frac12t$\\
\bottomrule
\end{tabular}
\caption{The slice measure of $T=T(\frac18,\frac25)$: the low zone on the left, the high zone on the right.}
\label{tab:aT}
\end{table}

\begin{proof}
(a) By~\eqref{eq:ahigh} with $\beta=\frac18$, $a_T(t)=\frac{233}{800}-\int_0^{t-3/8}a_E$ for
$t\ge\frac14$, where the integral is $0$ for $t\le\frac38$. By~\eqref{eq:alow}, for $t<\frac14$,
$a_T(t)=2\int_{\max(0,t-1/8)}^{t+1/8}r$, where $r=a_E$ on $[0,\frac14)$ and $r=\frac12$ on
$[\frac14,\frac12)$; note that $r=\frac12$ on all of $[\frac3{20},\frac12)$. Evaluating these integrals of
the piecewise linear function $a_E$ gives the table. As a check, $\int_0^{3/4}a_T=\frac{15803}{96000}$.

(b) We insert $b=a_T$, $f=\vol E(\frac5{12})=\frac7{24}$ and $\beta=\frac3{20}$ into~\eqref{eq:volplus}. The
three terms are the high zone $\int_0^{3/4}a_T(t)\min(\frac12,t+\frac1{10})\,dt$, the inherited rows
$2\int_0^{1/4}a_T(t)\,\omega_\ell(t)\,dt$ with $\omega_\ell(t)=t+\frac1{10}$ for $t<\frac3{20}$ and
$\frac25-t$ after, and the free rows $\frac9{400}\cdot\frac7{24}$.
\end{proof}

So $T_4$ gives the constant $(\frac{939039}{10240000})^{1/4}>0.55029$, by integrals of polynomials of degree at
most three that can be checked by hand. For comparison, the product of the set of~\cite{EPS} with itself has
volume $(7/24)^2\approx0.08507$. Other parameters are slightly better.

\begin{proposition}\label{prop:T4improved}
Let $\gamma=\frac{393}{1000}$ and $\beta_3=\beta_4=\frac{137}{1000}$. Then
$\vol T(\beta_3,\gamma)=\frac{16445487}{100000000}$, and the block
$T(\beta_3,\gamma)^+\bigl(E(\frac5{12}),\beta_4\bigr)$ in dimension~$4$ has volume
\[
  \frac{73446607289}{800000000000}\;=\;0.09180825911125 ,
\]
which gives the constant $0.550453>0.55045$.
\end{proposition}

\begin{proof}
The first value is $V_1(\frac{137}{1000},\frac{393}{1000})$. For the second we evaluated~\eqref{eq:ahigh},
\eqref{eq:alow} and~\eqref{eq:volplus} in exact integer arithmetic, with the programme described in the proof
of Proposition~\ref{prop:higher} on the grid of mesh $\frac1{1000}$. The set is a block by
Corollary~\ref{cor:T} and Theorem~\ref{thm:step}.
\end{proof}

We do not claim that these parameters are optimal.

\subsection*{Higher floors}

\begin{proposition}\label{prop:higher}
With the parameters of Definition~\ref{def:tower}, the volumes of the floors $T_5,\dots,T_{1000}$ are
rational numbers, of which we list the leading digits. The volumes and the constants $(\vol T_d)^{1/d}$ are
rounded down, and the ratios are rounded.
\begin{center}
\small
\renewcommand{\arraystretch}{1.2}
\begin{tabular}{@{}rlll@{}}
\toprule
$d$ & $\vol T_d$ & $(\vol T_d)^{1/d}$ & $\vol T_d/\vol T_{d-1}$\\
\midrule
$5$ & $5.08547\cdot10^{-2}$ & $0.551145$ & $0.554559$\\
$10$ & $2.67517\cdot10^{-3}$ & $0.553012$ & $0.554971$\\
$100$ & $2.58944\cdot10^{-26}$ & $0.554794$ & $0.554993$\\
$1000$ & $1.87077\cdot10^{-256}$ & $0.554973$ & $0.554993$\\
\bottomrule
\end{tabular}
\end{center}
In particular $\vol T_{1000}>0.55497321^{1000}$.
\end{proposition}

\begin{proof}
The recursion~\eqref{eq:ahigh}--\eqref{eq:alow} maps piecewise polynomial functions with rational
breakpoints and coefficients to functions of the same kind, so it can be evaluated exactly. All breakpoints
are multiples of $\frac1{120}$. On the cell $[\frac j{120},\frac{j+1}{120})$ we write the slice measure as
$a_{T_d}\bigl(\frac{j+u}{120}\bigr)=Q^{-1}\sum_iA_{j,i}\,u^i/i!$ with integers $A_{j,i}$ and one common
denominator $Q$. In this basis integration shifts the coefficients, and since $\beta_d$, $\frac14$ and
$\frac12$ are multiples of $\frac1{120}$, the limits of the integrals in~\eqref{eq:ahigh} and~\eqref{eq:alow}
move from cell to cell with the same local variable $u$. So one step consists of additions and subtractions of
integers, apart from the $90$ constants of integration. We carried this out with integers of arbitrary length
(programme \file{tower\_exact.c}). The exact values for $d=5$ and $d=6$ are listed in
Table~\ref{tab:constants}; the volume of $T_{1000}$ is a fraction whose numerator has $4088$ decimal digits, and
the last inequality is a comparison of integers.
\end{proof}

\subsection*{The growth rate of the tower}

For $d\ge4$ the step from $T_d$ to $T_{d+1}$ uses $\beta=\frac16$ and the free set $T_{d-1}$. By
Proposition~\ref{prop:volstep}, the pair $(a_{T_{d+1}},\vol T_d)$ is the image of the pair
$(a_{T_d},\vol T_{d-1})$ under the linear operator
\[
  \mathcal L(b,f)\;=\;\Bigl(Kb+f\,q,\ \int_0^{3/4}b(t)\,dt\Bigr),
\]
where
\[
  (Kb)(w)=\begin{cases}
  2\displaystyle\int_{\max(0,\,w-\frac1{12})}^{\min(\frac14,\,w+\frac16)}b(t)\,dt & \text{for }0\le w<\tfrac14,\\[2.5ex]
  \displaystyle\int_{\max(0,\,w-\frac13)}^{3/4}b(t)\,dt & \text{for }\tfrac14\le w<\tfrac34,
  \end{cases}
\]
and $q(w)=2\max(0,\,w-\tfrac1{12})$ for $w<\frac14$, $q(w)=0$ for $w\ge\frac14$.
The operator $\mathcal L$ is positive: if $b\ge0$ and $f\ge0$, then both components of $\mathcal L(b,f)$ are
non-negative. So inequalities between pairs, meant pointwise in both components, are preserved by
$\mathcal L$. This allows us to bound the growth of the volumes from both sides without computing them.

\begin{lemma}\label{lem:sandwich}
Let $u=(b,f)$ with a bounded measurable function $b\colon[0,\frac34)\to(0,\infty)$ that is bounded away from
$0$, and with $f>0$. Suppose that $\lambda_-u\le\mathcal Lu\le\lambda_+u$. Then there are constants
$c_1,c_2>0$ with
\[
  c_1\,\lambda_-^{\,d}\;\le\;\vol T_d\;\le\;c_2\,\lambda_+^{\,d}\qquad\text{for all }d\ge4 .
\]
\end{lemma}

\begin{proof}
For $d\ge3$ the slice measure $a_{T_d}$ is bounded, and it is bounded away from $0$ on $[0,\frac34)$. Indeed,
this holds for $a_E$, and if it holds for $b$ and if $f>0$, then by~\eqref{eq:ahigh} and~\eqref{eq:alow} the
slice measure of $B^+$ at every level is the integral of a function that is bounded away from $0$ over an
interval of length at least $\min(\beta,\frac14)>0$. So there are $c_1,c_2>0$ with
$c_1\lambda_-^4u\le(a_{T_5},\vol T_4)\le c_2\lambda_+^4u$. Applying $\mathcal L$ repeatedly and using its
positivity gives $c_1\lambda_-^{\,d}u\le(a_{T_{d+1}},\vol T_d)\le c_2\lambda_+^{\,d}u$ for $d\ge4$, and the second
component, divided by $f$, is the claim (with other constants).
\end{proof}

\begin{proposition}\label{prop:cert}
There is a pair $u$ as in Lemma~\ref{lem:sandwich} with
\[
  0.55499308\,u\;\le\;\mathcal Lu\;\le\;0.55499309\,u .
\]
\end{proposition}

\begin{proof}
Let $N=48000$. We take $u=(b,f)$ with a function $b$ that is linear on each of the $36000$ cells
$[\frac iN,\frac{i+1}N)$, with positive integer values at the ends of the cells and with separate values on the
two sides of $\frac14$, and with a positive integer $f$. Since $\frac1{12}$, $\frac16$, $\frac14$ and $\frac13$
are multiples of $\frac1N$, the function $\mathcal Lu$ is a quadratic polynomial on each cell, and so are
$\mathcal Lu-\lambda_-u$ and $\lambda_+u-\mathcal Lu$. A quadratic polynomial $p$ is non-negative on $[0,1]$ if
its three Bernstein coefficients $p(0)$, $2p(\frac12)-\frac12(p(0)+p(1))$ and $p(1)$ are non-negative. For
suitable values of $u$, found by a numerical iteration of $\mathcal L$, these coefficients are non-negative
integers, up to a common factor, on every cell, and the second component of $\mathcal Lu$, which is $\int b$,
satisfies both inequalities as well (programme \file{spectral\_cert.c}, in integer arithmetic).
\end{proof}

\begin{proof}[Proof of Theorem~\ref{thm:main}]
The sets $T_d$ are reducible by Theorem~\ref{thm:tower}. By Lemma~\ref{lem:sandwich} and
Proposition~\ref{prop:cert},
\[
  0.55499308\;\le\;\liminf_{d\to\infty}(\vol T_d)^{1/d}\;\le\;\limsup_{d\to\infty}(\vol T_d)^{1/d}\;\le\;0.55499309 .
\]
It remains to show that the ratios $\vol T_{d+1}/\vol T_d$ converge; the constants $(\vol T_d)^{1/d}$ then
converge to the same limit. The slice measures $a_{T_d}$, $d\ge3$, are continuous on $[0,\frac14)$ and on
$[\frac14,\frac34)$ and have limits at the right ends. So we may consider $\mathcal L$ on the Banach space
$X=C([0,\frac14])\times C([\frac14,\frac34])\times\mathbb{R}$, ordered by the cone $X_+$ of triples of
non-negative functions and numbers; the interior of $X_+$ consists of the triples that are positive
everywhere. The operator $\mathcal L$ maps bounded sets to sets of uniformly bounded functions with a common
Lipschitz constant, so it is compact by the theorem of Arzel\`a and Ascoli.

We claim that $\mathcal L^4$ maps every non-zero element of $X_+$ into the interior of $X_+$. Let
$x=(b,f)\in X_+$ with $b\ne0$. Then $\int b>0$, so $\mathcal Lx$ has a positive last component, and it is
positive at the levels $w\in[\frac14,\frac13]$, where the lower limit of the integral is $0$. In
$\mathcal L^2x$ the last component is positive, the levels $w\in(\frac1{12},\frac14]$ are positive because of
the term $f\,q$, and the levels $w\in[\frac14,\frac23)$ are positive, because there the interval
$[\max(0,w-\frac13),\frac34]$ meets $[\frac14,\frac13]$ in an interval. In $\mathcal L^3x$ all levels below
$\frac14$ are positive, since every window $[\max(0,w-\frac1{12}),\min(\frac14,w+\frac16)]$ meets
$(\frac1{12},\frac14]$ in an interval, and all levels $w\in[\frac14,\frac34]$ are positive, since
$w-\frac13\le\frac5{12}<\frac23$. So $\mathcal L^3x$ lies in the interior of $X_+$, and so does
$\mathcal L^4x$. If $b=0$ and $f>0$, then $\mathcal Lx=(f\,q,0)$ has a non-zero first component, and the
previous argument applies to $\mathcal Lx$.

By the theorem of Krein and Rutman~\cite[Theorem~6.3]{KreinRutman}, the compact operator $\mathcal L^4$ has
an algebraically simple eigenvalue $\rho>0$ with an eigenvector $u_\infty$ in the interior of $X_+$, and all its other
eigenvalues have smaller absolute value. Since $\mathcal L$ commutes with $\mathcal L^4$, the vector
$\mathcal Lu_\infty$ is again an eigenvector of $\mathcal L^4$ for $\rho$, so
$\mathcal Lu_\infty=r_\infty u_\infty$ with $r_\infty>0$ and $r_\infty^4=\rho$. If $\mu$ is an eigenvalue of
$\mathcal L$ with $|\mu|=r_\infty$, then $\mu^4$ is an eigenvalue of $\mathcal L^4$ of absolute value $\rho$,
so $\mu^4=\rho$, its eigenvectors are multiples of $u_\infty$, and $\mu=r_\infty$. Hence $r_\infty$ is an algebraically simple
eigenvalue of $\mathcal L$ (a generalised eigenvector of $\mathcal L$ is one of $\mathcal L^4$), and the rest of the spectrum of the compact operator $\mathcal L$ lies in a disc
of radius smaller than $r_\infty$. It follows that $r_\infty^{-k}\mathcal L^kx$ converges to a multiple
$\ell(x)\,u_\infty$ for every $x\in X$, and $\ell(x)>0$ for non-zero $x\in X_+$, because
$\mathcal L^4x\ge\varepsilon u_\infty$ for some $\varepsilon>0$. Applied to $x=(a_{T_5},\vol T_4)$ this gives
$\vol T_{d+1}/\vol T_d\to r_\infty$, which proves~(a). Hence there is a dimension $d$ with
$(\vol T_d)^{1/d}>0.554993$, and Corollary~\ref{cor:torus} gives~(b). Part~(c) is contained in
Proposition~\ref{prop:higher}, again with Corollary~\ref{cor:torus}.
\end{proof}

The proof of part~(b) does not name the dimension $d$. The ratios converge quickly, but the constants
$(\vol T_d)^{1/d}$ converge slowly, as the table shows, so $d$ has to be large.

\begin{remark}\label{rem:limit}
A change of finitely many parameters or free sets does not change the limit $r_\infty$, since the state of the
recursion after these floors is again a pair of a positive bounded function and a positive number. If another
parameter $\beta$ is used on all high floors, the growth rate changes. A floating-point evaluation of the
recursion gives the following approximate limits of the ratios.
\begin{center}
\small
\renewcommand{\arraystretch}{1.2}
\begin{tabular}{@{}lllllllll@{}}
\toprule
$\beta$ & $0.12$ & $0.14$ & $0.15$ & $0.16$ & $1/6$ & $0.17$ & $0.18$ & $0.20$\\
ratio & $0.53806$ & $0.54943$ & $0.55290$ & $0.55470$ & $0.55499$ & $0.55487$ & $0.55354$ &
$0.54702$\\
\bottomrule
\end{tabular}
\end{center}
Our proof of convergence and our enclosure concern $\beta=\frac16$ only; the other values are numerical, and
we have not optimised over sequences of parameters.
\end{remark}

\section{Larger sets in dimensions two, three and four}\label{sec:dim3}

In this section $0<\beta\le\frac14$ and $\frac14\le\gamma\le\frac12$. We treat two sets side by side: in
dimension~$3$ the second floor $T=T(\beta,\gamma)$, and in dimension~$2$ its analogue one dimension lower, the
set $[0,\frac12)^+$ of Remark~\ref{rem:floor1}. For $\frac12\le t<\frac34$ the high zones of both sets are cut
off below the line $s=t-\frac12+\beta$, and a notch remains free. Part of this notch can be added. The price
is a removal order in which the time of removal depends on $s$ as well as on $t$, so the resulting sets are
not blocks. The arguments are the same in both dimensions; the additional coordinate $u$ of dimension~$3$
enters only at two places.

\begin{definition}\label{def:Mstar}
Let
\[
  \phi(\tau)=\max\bigl(0,\tau-\tfrac14+\beta\bigr),\qquad g(\tau)=\frac{(\beta-\tau)(1-8\tau)}{1+8\tau},
\]
and let $N\subseteq[0,1)^2$ be the set of points $(s,t)$ with $\frac12\le t<\frac12+\min(\beta,\frac18)$ and
$\phi(\tau)\le s<g(\tau)$, where $\tau=t-\frac12$.
\begin{enumerate}
\item[(a)] Let $R_\beta\subseteq[0,1)^2$, with coordinates $(s,t)$ and level $t$, be the union of $N$, the
\emph{low zone} $L$ and the \emph{high zone} $H$, where
\begin{align*}
  L&=\bigl\{(s,t):\ 0\le t<\tfrac14,\ \ \max(0,\,t-\tfrac14+\beta)\le\hf{s}<t+\beta\bigr\},\\
  H&=\bigl\{(s,t):\ \tfrac14\le t<\tfrac34,\ \ \max(0,\,t-\tfrac12+\beta)\le s<\tfrac12\bigr\}.
\end{align*}
\item[(b)] Let $M^*$ be the set of points $(u,s,t)\in[0,1)^3$ with $(s,t)\in N$ and $\hf{u}<\rho(s)$.
\end{enumerate}
We call $N$ the \emph{notch} of $R_\beta$ and $M^*$ the notch of $T\cup M^*$.
\end{definition}

\begin{figure}[tbp]
\centering
\includegraphics[width=0.5\textwidth]{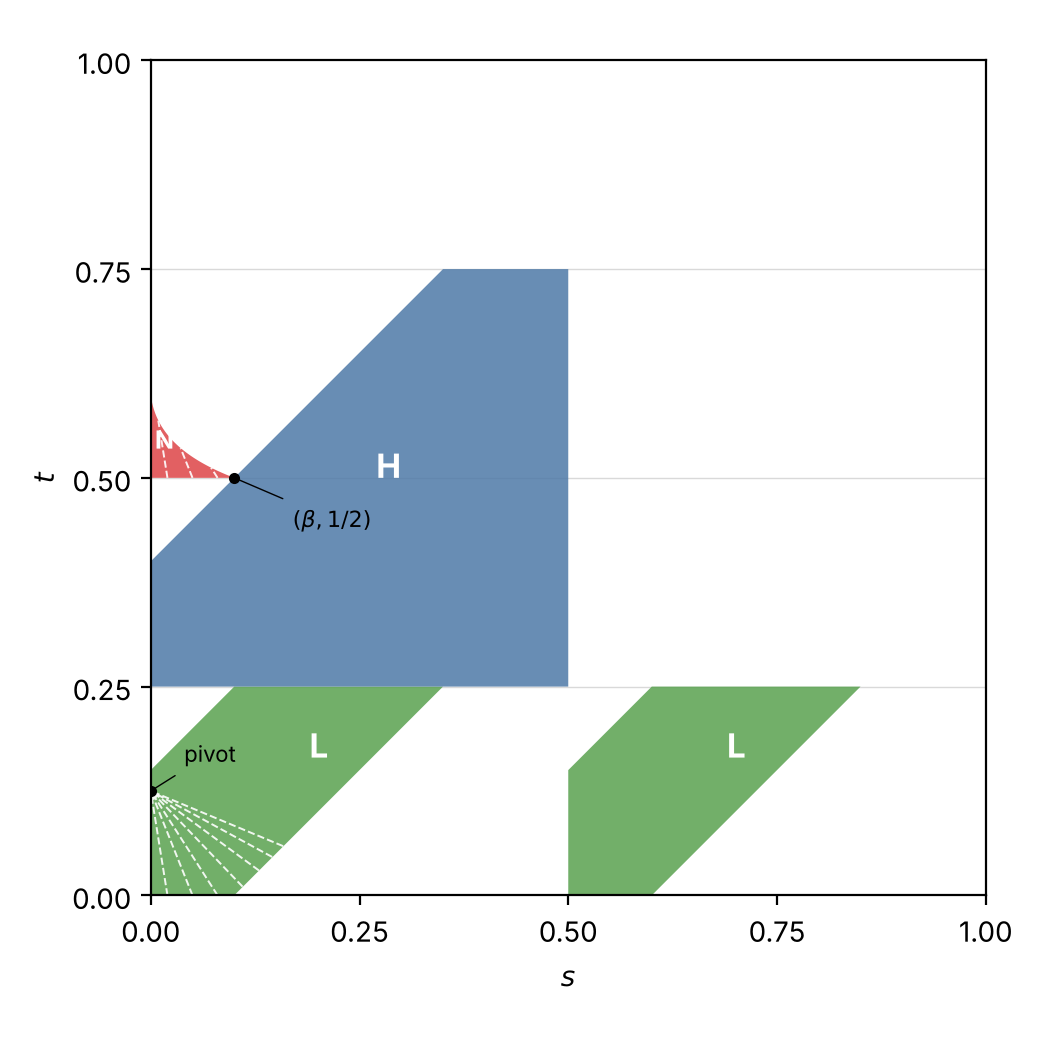}
\caption{The set $R_\beta$ for $\beta=\frac1{10}$. The low zone $L$ (green) consists of two copies, the high
zone $H$ (blue) is cut off below the line $s=t-\frac12+\beta$, and the notch $N$ (red) lies below this cut.
The dashed lines are lines of the fan through the pivot (Definition~\ref{def:fan}).}
\label{fig:R}
\end{figure}

\begin{figure}[tbp]
\centering
\includegraphics[width=\textwidth]{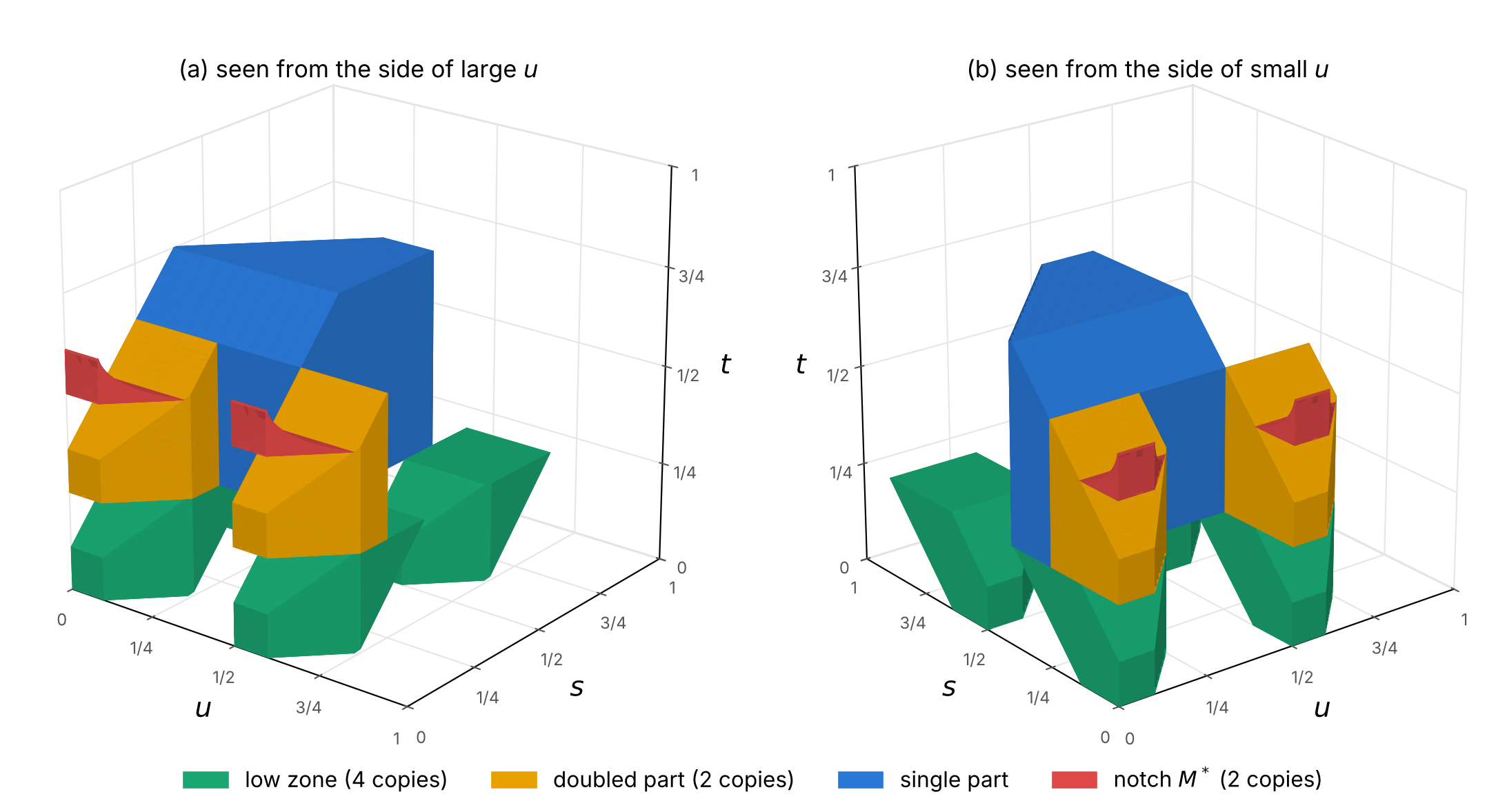}
\caption{The set $T(\frac7{50},\frac25)\cup M^*$ of Theorem~\ref{thm:dim3}(b) in the unit cube with coordinates
$(u,s,t)$, seen from two sides. The low zone (green) consists of four copies, the doubled part (yellow)
and the notch $M^*$ (red) of two each, and the single part (blue) is one piece. In (a) the widths in $u$ show
the profile $\rho$; in (b) the faces towards small $u$ show the shapes in the $(s,t)$-plane, among them the cut
of the high zone and the curved boundary of the notch.}
\label{fig:structure}
\end{figure}

For $\beta\le\frac18$ we have $\phi(\tau)=0$ for all $\tau<\beta$, and the notch is the region $0\le s<g(\tau)$,
$0\le\tau<\beta$. For $\beta>\frac18$ the notch is truncated: the curves $\phi$ and $g$ meet at
$\tau=1/(64\beta)<\frac18$, and the rows above this level are empty. The set $L\cup H$ is the set
$[0,\frac12)^+$ of Remark~\ref{rem:floor1}, so it is a block by Theorem~\ref{thm:step}, and it has the slice
measures and the area of $E(\frac12-\beta)$. On $N$ we have $s<g(\tau)\le\beta-\tau\le t-\frac12+\beta$, so the
notch lies below the cut of the high zone. Hence $N$ is disjoint from $L\cup H$, and $M^*$ is disjoint from
$T$. Figure~\ref{fig:R} shows $R_{1/10}$, and Figure~\ref{fig:structure} shows $T\cup M^*$ for
$\beta=\frac7{50}$ and $\gamma=\frac25$, the set of Theorem~\ref{thm:dim3}(b).

\subsection*{Chart coordinates and the fan order}
For a point of the low zone or of the notch we write $\sigma=\hf{s}$ and $\tau=\hf{t}$. So $\tau=t$ on the
low zone and $\tau=t-\frac12$ on the notch, and $\sigma=s$ on the notch. In both dimensions the points
$(\sigma,\tau)$ lie in the region $\Omega=\Omega_\beta$ of Lemma~\ref{lem:chart}; for the notch this holds
because $\tau<\frac18$ and $\phi(\tau)\le\sigma<\tau+\beta$. In these coordinates the low zone consists of two
copies of $\Omega$, which differ by $\frac12$ in $s$, and the notch is part of a third copy, which differs from
them by $\frac12$ in $t$. In dimension~$3$ a point of any of these copies has the fibre
$\hf{u}<\rho(\sigma)$.

\begin{definition}\label{def:fan}
The \emph{fan order} on $R_\beta$ and on $T\cup M^*$ removes the points in three stages.
\begin{itemize}[leftmargin=4.8em]
\item[Stage 0:] all points of the low zone and of the notch with $\tau<\frac18$, by increasing
$\theta=\sigma/(1-8\tau)$; ties are broken by smaller $\tau$ first, then, in dimension~$3$, by smaller
$\hf{u}$ first.
\item[Stage 1:] the remaining points of the low zone, by increasing $\tau$, then by increasing $\sigma$,
then, in dimension~$3$, by increasing $\hf{u}$.
\item[Stage 2:] the high zone; in dimension~$2$ by increasing $t$ and then by decreasing $s$, in
dimension~$3$ in the removal order of Corollary~\ref{cor:T}.
\end{itemize}
Ties that remain, between points of different copies with the same chart coordinates, are broken arbitrarily.
\end{definition}

\begin{figure}[tbp]
\centering
\includegraphics[width=\textwidth]{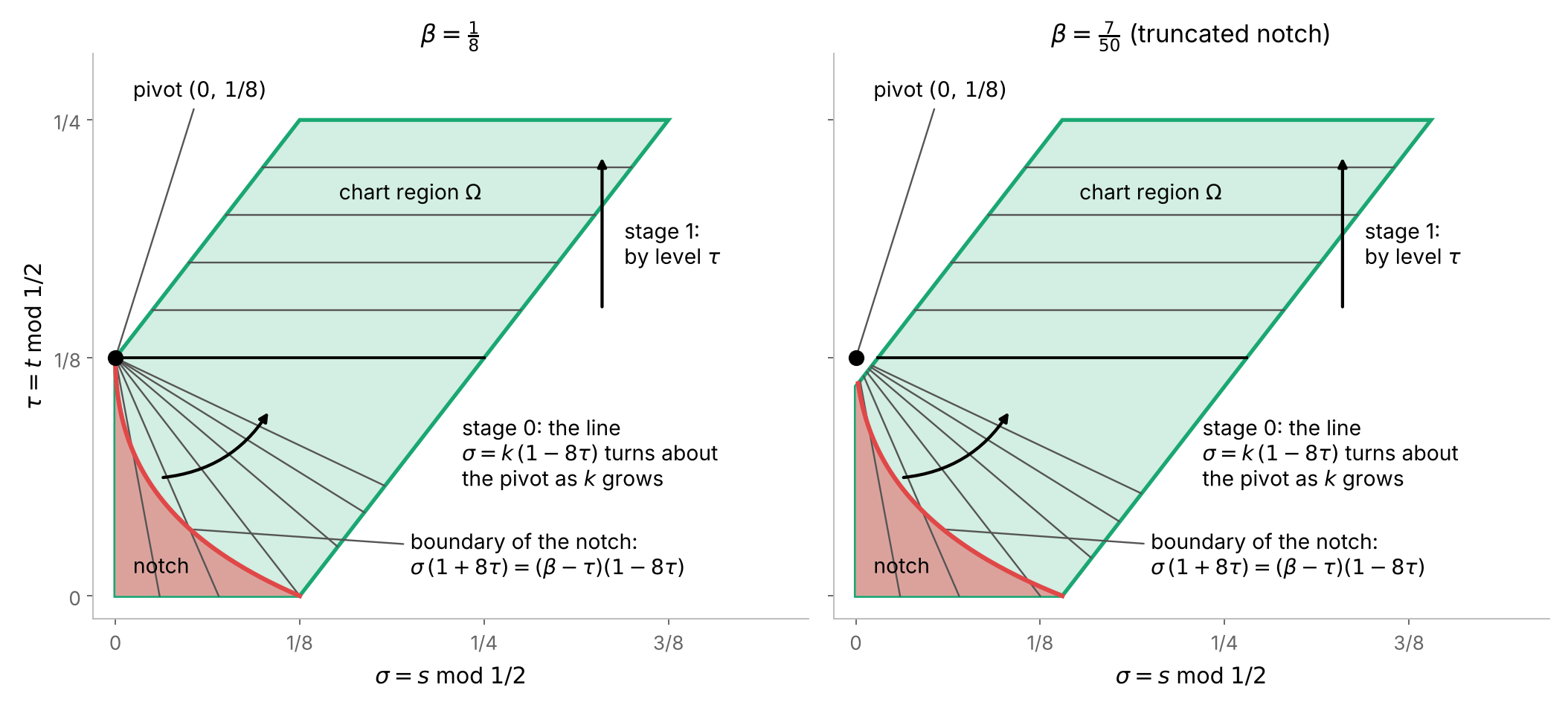}
\caption{The chart region $\Omega$ with the notch (red) and the fan order, for $\beta=\frac18$ and for
$\beta=\frac7{50}$. In stage~0 the line $\sigma=k(1-8\tau)$ turns about the pivot; in stage~1 the remaining
levels are removed from below. For $\beta>\frac18$ the lower end $\sigma=\tau-\frac14+\beta$ of the window is
positive near the pivot, so the pivot lies outside $\Omega$, and the notch is truncated where its
curved boundary meets this line, at $\tau=1/(64\beta)$.}
\label{fig:fan}
\end{figure}

The set $\theta=k$ is the line $\sigma=k(1-8\tau)$ through the \emph{pivot} $(\sigma,\tau)=(0,\frac18)$. As
$k$ grows, the line turns about the pivot and sweeps the part of $\Omega$ below $\tau=\frac18$
(Figure~\ref{fig:fan}). Every point of the notch is removed in stage~0.

\begin{theorem}\label{thm:Mstar}
Let $0<\beta\le\frac14$ and $\frac14\le\gamma\le\frac12$. Then $R_\beta$ and $T(\beta,\gamma)\cup M^*$ are
non-mid-point reducible. More precisely, if a non-empty finite subset $S$ of one of these sets meets the low
zone or the notch, then every point of $S$ that comes first in the fan order is a non-mid-point of $S$.
\end{theorem}

\begin{proof}
We give the proof for $T\cup M^*$ and point out the two places where the coordinate $u$ is used; for $R_\beta$
the same argument applies with $u$ omitted. If $S$ lies in the high zone, it is a subset of the block $L\cup H$
or $T$, so it has a non-mid-point. Otherwise let $y$ be a point of $S$ that comes
first in the fan order; it lies in the low zone or in the notch. Suppose that $2y\equiv x+z\pmod 1$ with
distinct $x,z\in S$. Then neither $x$ nor $z$ comes before $y$. We distinguish three cases.

\emph{Case 1: $x$ and $z$ lie in the high zone.} Then $t_x+t_z\in[\frac12,\frac32)$ and $2t_y\equiv2\tau_y$
with $2\tau_y\in[0,\frac12)$, so $t_x+t_z=2\tau_y+1$. On the low zone and on the notch we have
$\sigma_y<\tau_y+\beta$. Together with the cuts $s_x\ge t_x-\frac12+\beta$ and $s_z\ge t_z-\frac12+\beta$ this
gives
\[
  2\sigma_y\;<\;2\tau_y+2\beta\;=\;(t_x-\tfrac12+\beta)+(t_z-\tfrac12+\beta)\;\le\;s_x+s_z .
\]
For $R_\beta$ we have $s_x,s_z<\frac12$, so $2\sigma_y$ and $s_x+s_z$ are different numbers in $[0,1)$, and
since $2s_y\equiv2\sigma_y$ the congruence $2s_y\equiv s_x+s_z$ fails. For $T\cup M^*$ we have
$s_x+s_z<\frac32$, so the congruence would give $s_x+s_z=2\sigma_y+1$. Then $s_x,s_z>\frac14$, so $x$ and $z$
lie in the single part, and, since $\hf{u_y}<\rho(\sigma_y)\le\sigma_y+\frac12-\gamma$,
\[
  u_x+u_z\;\ge\;s_x+s_z-2\gamma\;=\;2\sigma_y+1-2\gamma\;>\;2\hf{u_y},
\]
while $u_x+u_z<1$. So the congruence fails in $s$ or in $u$, a contradiction.

\emph{Case 2: $x$ and $z$ lie in the low zone or in the notch.} By Lemma~\ref{lem:chart} the chart point of
$y$ is the mid-point of the chart points of $x$ and $z$, and $x$ and $z$ lie in the same copy of $\Omega$.
Suppose first that the chart points of $x$ and $z$ are different.
\begin{itemize}
\item If $\tau_y\ge\frac18$, then $y$ is removed in stage~1. If $\tau_x\ne\tau_z$, the point with the
smaller $\tau$ lies below level $\tau_y$ and comes before $y$. If $\tau_x=\tau_z$, then both equal
$\tau_y$ and $\sigma_x\ne\sigma_z$, and the point with the smaller $\sigma$ comes before $y$.
\item If $\tau_y<\frac18$, put $k=\theta(y)$ and $\ell(\sigma,\tau)=\sigma-k(1-8\tau)$. This function is affine
and vanishes at the chart point of $y$, so $\ell(x)+\ell(z)=0$. A chart point with $\tau\ge\frac18$ has
$\ell\ge0$, and for a chart point with $\tau<\frac18$ the sign of $\ell$ is the sign of $\theta-k$. If
$\ell(x)<0$, then $\tau_x<\frac18$ and $\theta(x)<k$, so $x$ comes before $y$; likewise for $z$. If
$\ell(x)=\ell(z)=0$, then $\tau_x\ne\tau_z$, since two different points of the line $\ell=0$ have different
$\tau$. The point with the smaller $\tau$, say $x$, has $\tau_x<\tau_y<\frac18$ and $\theta(x)=k$, and it
comes before $y$ by the tie rule.
\end{itemize}
Either way one of $x,z$ comes before $y$, a contradiction. Now suppose that $x$ and $z$ have the same chart
point; it is then the chart point of $y$, and $x$ and $z$ have the same $s$ and $t$. For $R_\beta$ this means
$x=z$, a contradiction. For $T\cup M^*$ the points $x$ and $z$ differ in $u$. By Lemma~\ref{lem:interval},
applied to the $u$-coordinates with $I=[0,\frac14)$, the numbers $u_x$ and $u_z$ lie in the same half of
$[0,1)$ and $\hf{u_y}$ is the average of $\hf{u_x}$ and $\hf{u_z}$. These two values are different, so the
point with the smaller one comes before $y$, again a contradiction.

\emph{Case 3: one of $x,z$, say $z$, lies in the high zone, and the other in the low zone or in the notch.}
Suppose that $x$ lies in the low zone. Then $t_z\equiv2\tau_y-\tau_x\pmod 1$, and $2\tau_y-\tau_x$ lies in
$(-\frac14,\frac12)$, so $t_z=2\tau_y-\tau_x\ge\frac14$. Hence $\tau_y\ge\frac18+\frac12\tau_x$ and
$\tau_x\le2\tau_y-\frac14<\tau_y$. So $y$ is removed in stage~1, and $x$ lies at a lower level and comes
before $y$, a contradiction.

Suppose that $x$ lies in the notch. Then $x$ is removed in stage~0, hence so is $y$, and $\theta(y)\le k$ for
$k=\theta(x)$. So $\tau_y<\frac18$ and $\sigma_y\le k(1-8\tau_y)$, while $\sigma_x=k(1-8\tau_x)$. Now
$t_z\equiv2\tau_y-\tau_x+\frac12\pmod 1$, and $2\tau_y-\tau_x$ lies in $(-\frac18,\frac14)$, so
$t_z=2\tau_y-\tau_x+\frac12$. Also $s_z\equiv2\sigma_y-\sigma_x$. Since $\sigma_x<\beta\le\frac14$ and
$\sigma_y<\tau_y+\beta<\frac38$, the number $2\sigma_y-\sigma_x$ lies in $(-\frac14,\frac34)$, and a negative
value would put $s_z$ above $\frac34$, outside the high zone. So $s_z=2\sigma_y-\sigma_x$, and
the cut $s_z\ge t_z-\frac12+\beta$ of the high zone reads $2\sigma_y-\sigma_x\ge2\tau_y-\tau_x+\beta$.
Inserting the expressions for $\sigma_y$ and $\sigma_x$ gives
\[
  k\,(1+8\tau_x)\;\ge\;\beta-\tau_x+2\tau_y\,(1+8k)\;\ge\;\beta-\tau_x .
\]
But $\sigma_x<g(\tau_x)$, divided by $1-8\tau_x>0$, says $k(1+8\tau_x)<\beta-\tau_x$, a contradiction.
\end{proof}

\begin{proposition}\label{prop:volM}
\begin{enumerate}
\item[(a)] For $\beta\le\frac18$ the notch $N$ has area
$\int_0^\beta g(\tau)\,d\tau=\frac{1+8\beta}{32}\ln(1+8\beta)-\frac\beta4-\frac{\beta^2}2$, and
\[
  \vol R_\beta\;=\;\frac9{32}-2\beta^2+\frac{1+8\beta}{32}\,\ln(1+8\beta).
\]
For $\beta=\frac1{10}$ this is $\frac{209}{800}+\frac9{160}\ln\frac95\approx0.2943130$. On $(0,\frac18]$ the
area is strictly concave. It is largest at $\beta^*=\frac18(x^*-1)\approx0.098942$, where
$x^*\approx1.791537$ is the solution $x>1$ of $\ln x=2x-3$, and its largest value is
$\frac1{32}\bigl((x^*)^2-x^*+8\bigr)\approx0.2943146$.
\item[(b)] Let $\beta\le\gamma-\frac14$, and let $\tau^*=\beta$ if $\beta\le\frac18$ and $\tau^*=1/(64\beta)$
otherwise. Then
\[
  \vol M^*\;=\;\int_0^{\tau^*}\bigl(g(\tau)^2+(1-2\gamma)\,g(\tau)\bigr)\,d\tau
  -\int_0^{\tau^*}\bigl(\phi(\tau)^2+(1-2\gamma)\,\phi(\tau)\bigr)\,d\tau .
\]
For $\gamma=\frac25$ this equals $\frac{55}{1536}-\frac{\ln 2}{20}\approx0.001150$ if $\beta=\frac18$, and
$\frac{12254387}{336000000}-\frac{4399}{80000}\ln\frac{53}{28}\approx0.0013846$ if $\beta=\frac7{50}$.
\end{enumerate}
\end{proposition}

\begin{proof}
(a) The set $L\cup H$ has the area $\vol E(\frac12-\beta)=\frac9{32}+\frac\beta4-\frac32\beta^2$. Since
$\frac{1-8\tau}{1+8\tau}=\frac2{1+8\tau}-1$ and $2(\beta-\tau)=\frac{1+8\beta}4-\frac{1+8\tau}4$, we have
\[
  g(\tau)\;=\;\frac{1+8\beta}{4(1+8\tau)}-\frac14-(\beta-\tau),
\]
and integration gives the area of $N$ and the formula for $A(\beta)=\vol R_\beta$. We have
$A'(\beta)=-4\beta+\frac14\ln(1+8\beta)+\frac14$ and $A''(\beta)=-4+\frac2{1+8\beta}<0$. So $A'$ vanishes
exactly where $\ln x=2x-3$ with $x=1+8\beta$, and inserting this gives
$A=\frac1{32}\bigl(9-(x-1)^2+x(2x-3)\bigr)=\frac1{32}(x^2-x+8)$.

(b) On $M^*$ we have $s<\beta\le\gamma-\frac14$, so $\rho(s)=s+\frac12-\gamma$ and the fibre has measure
$2(s+\frac12-\gamma)$. Integrating over $\phi(\tau)\le s<g(\tau)$ gives the two integrands; for
$\beta>\frac18$ the curves $\phi$ and $g$ meet where $16\beta\tau=\frac14$. The values follow with the
decomposition of $g$ from~(a), with $1+8\tau^*=2$ and $1+8\tau^*=\frac{53}{28}$, respectively.
\end{proof}

\begin{proof}[Proof of Theorem~\ref{thm:dim3}(a) and (b)]
(a) By Theorem~\ref{thm:Mstar} the set $R_{1/10}$ is reducible, and by Proposition~\ref{prop:volM}(a) its area
is $v=\frac{209}{800}+\frac9{160}\ln\frac95$. Corollary~\ref{cor:torus} gives the constant $v^{1/2}>0.54250$.

(b) By Theorem~\ref{thm:Mstar} the set $T(\frac7{50},\frac25)\cup M^*$ is reducible. By
Propositions~\ref{prop:volT} and~\ref{prop:volM}(b) its volume is
\[
  V_1(\tfrac7{50},\tfrac25)+\vol M^*=\tfrac{7711}{46875}+\tfrac{12254387}{336000000}-\tfrac{4399}{80000}\ln\tfrac{53}{28}
  =\tfrac{4501789}{22400000}-\tfrac{4399}{80000}\ln\tfrac{53}{28}=0.1658858\ldots,
\]
with cube root $0.549460\ldots>0.54946$. The constants $C$ are obtained in Section~\ref{sec:behrend}.
\end{proof}

Here and below, inequalities that involve logarithms are verified with the series
$\ln x=2\sum_{j\ge0}z^{2j+1}/(2j+1)$, $z=\frac{x-1}{x+1}$, whose partial sums are lower bounds and whose
remainder is bounded by a geometric series; for example $0.58773<\ln\frac95<0.58779$ from three terms with
$z=\frac27$.

For $\beta=\frac18$ and $\gamma=\frac25$, the parameters of the tower, the set $T_3\cup M^*$ has the volume
$\frac{12827}{64000}-\frac{\ln 2}{20}\approx0.165765$ and gives $c>0.54932$. A numerical optimisation over $\beta$
and $\gamma$ gives the volume $\approx0.165888$ near $(\beta,\gamma)=(0.138,0.398)$; we do not claim that this is
optimal. In dimension~$2$ the truncated notch gives nothing new: numerically the area of $R_\beta$ decreases for
$\beta>\frac18$.

\begin{remark}\label{rem:notblock}
Neither $R_\beta$ nor $T\cup M^*$ is a block. The points $x=(\frac\beta2,\frac12)\in N$,
$y=(\frac{3\beta}4,0)\in L$ and $z=(\beta,\frac12)\in H$ satisfy $t_x+t_z=2t_y+1$ and $x+z\equiv2y\pmod 1$, and
in dimension~$3$ the same holds after adding the coordinate $u=0$ to each point. So property~(B) fails, and
neither set can take the place of $B$ in the induction step. Both can take the place of the free set $F$,
which only needs to be reducible.
\end{remark}

\subsection*{The fourth floor with a new free set}
Let $\widetilde T_4=T\bigl(\frac{137}{1000},\frac{393}{1000}\bigr)^+\bigl(R_{1/10},\frac{137}{1000}\bigr)$. Its
free rows are copies of $R_{1/10}$, with the coordinates $(u,s)$ in the roles of $(s,t)$.

\begin{proof}[Proof of Theorem~\ref{thm:dim4}(c)]
The set $T(\frac{137}{1000},\frac{393}{1000})$ is a block by Corollary~\ref{cor:T}, and $R_{1/10}$ is reducible by
Theorem~\ref{thm:Mstar}. So $\widetilde T_4$ is a block by Theorem~\ref{thm:step}. By~\eqref{eq:volplus}, the
free set enters the volume of $B^+(F,\beta)$ only through the term $\beta^2\vol F$. So by
Proposition~\ref{prop:T4improved}
\begin{align*}
  \vol\widetilde T_4&=\frac{73446607289}{800000000000}+\frac{18769}{1000000}\Bigl(\frac9{160}\ln\frac95-\frac{73}{2400}\Bigr)\\
  &=\frac{218969684867}{2400000000000}+\frac{168921}{160000000}\ln\frac95\;\approx\;0.0918579 ,
\end{align*}
and its fourth root is $0.550527\ldots>0.55052$.
\end{proof}

With the parameters of the tower, the set $T_3^+(R_{1/10},\frac3{20})$ has the volume
$\frac{932031}{10240000}+\frac{81}{64000}\ln\frac95\approx0.0917626$, by Proposition~\ref{prop:volT4}, and gives
$c>0.55038$; this value can be checked by hand.

\begin{remark}\label{rem:sharp}
Case~3 is the only place in the proof of Theorem~\ref{thm:Mstar} where the shape of the notch enters, and the
sufficient inequality that defines $g$ is tight there: for a point $x$ on the curve, the critical point $y$
lies at level $\tau=0$ on the same line of the fan. This does not show that the notch is the largest region
that can be added, since the inequalities in the proof are sufficient for the fan order but not necessary for
other removal orders.
\end{remark}

The notch cannot simply be added to the block of the next floor. For a set $B$ that is not a block we define
$B^+(F,\beta)$ by the formulas of Definition~\ref{def:step}.

\begin{proposition}\label{prop:Mnext}
Let $\beta=\frac18$ and $\gamma=\frac25$, and let $x_0$ be one of the points
$(\frac1{20},\frac1{100},\frac{11}{20})$ and $(\frac1{10},\frac1{10},\frac{101}{200})$ of $M^*$. Then the set
$(T_3\cup\{x_0\})^+\bigl(E(\frac5{12}),\frac3{20}\bigr)$, which is the fourth floor $T_4$ together with the lifts
of $x_0$ to the levels of its high zone, is not reducible.
\end{proposition}

\begin{proof}
Let $m=13$ and choose the offset $\alpha\in[0,1)^4$ with $\alpha\equiv x_0\pmod{\frac1{13}}$ in the first three
coordinates and $\alpha_4=\frac3{10}$, so that $x_0$ and its lifts lie on the grid. Removing non-mid-points
from the discretisation $S_\alpha$ of the set, which has $2437$ points for the first and $2627$ for the second
point, ends with a set of $2208$, respectively $2308$, points, each of which is the mid-point of two others. So
$S_\alpha$ is not reducible, and by Lemma~\ref{lem:transfer} neither is the set.
\end{proof}

The proposition concerns these two points and the formulas of Definition~\ref{def:step} with the stated
parameters and free set. It does not exclude that a different induction mechanism could use $T_3\cup M^*$.
Heuristically, the obstacle is property~(B): the triples of Remark~\ref{rem:notblock} lift to triples of the
next floor whose levels wrap around once, with the middle point in the low zone. Within one floor such triples
only force an order, $x$ before $y$; in the next floor they stand in the way of removing the low zone before
the high zone. So the notches enter the higher floors only through the free sets.

\section{Potentials and Behrend's construction}\label{sec:behrend}

Elsholtz, Hunter, Proske and Sauermann~\cite{EHPS} improved Behrend's construction with two-dimensional
building blocks, sets that carry a function with the following property. In this section we show that the induction step of
Section~\ref{sec:step} carries such functions from one floor to the next, at the price of removing parts of
arbitrarily small volume, and we prove Theorem~\ref{thm:behrend}.

\begin{definition}\label{def:potential}
Let $R\subseteq[0,1)^d$ be measurable. A \emph{potential} on $R$ is a bounded measurable function
$f\colon R\to[0,\infty)$ such that
\[
  f(x)+f(z)\;\ge\;2f(y)+|x-z|^2\qquad\text{whenever }x,y,z\in R\text{ and }x+z\equiv2y\pmod 1 .
\]
Here $|x-z|$ is the Euclidean norm of the difference of $x$ and $z$, taken in $\mathbb{R}^d$. We say that $R$
\emph{carries} a potential if there is one.
\end{definition}

This is the condition of~\cite[Proposition~2.2]{EHPS}. The points $x,y,z$ need not be distinct. For $x=z$
the condition says that $f(x)\ge f(y)$ whenever $2y\equiv2x$, so points that differ by a vector in
$(\frac12\Z)^d$, the analogues of the cousins of~\cite[Definition~2.2]{EPS}, have the same value. A set that carries a potential is reducible: in a non-empty finite
subset, a point $y$ at which $f$ is largest is a non-mid-point, because $x+z\equiv2y$ with $x\ne z$ would
give $f(x)+f(z)>2f(y)$. So a potential is a quantitative form of a removal order, in which points with
larger values are removed earlier. We shall often use that for real numbers $a,c$ and $b=\frac12(a+c)$
\begin{equation}\label{eq:midquad}
  \kappa a^2+\kappa c^2-2\kappa b^2\;=\;\tfrac{\kappa}2\,(a-c)^2 .
\end{equation}
For example, by the second part of Lemma~\ref{lem:interval} and~\eqref{eq:midquad}, the function $2x^2$ is a
potential on the interval $[0,\frac12)$; this is the potential behind Behrend's construction.

\begin{proposition}[{cf.~\cite[Proposition~2.1 and Section~3]{EHPS}}]\label{prop:behrend}
Let $R\subseteq[0,1)^d$ be a measurable set of positive volume that carries a potential. Then
\[
  r_3(N)\;\ge\;N\cdot2^{-(C_R+o(1))\sqrt{\log_2N}}\qquad\text{with}\qquad
  C_R=2\sqrt{\tfrac2d\log_2\bigl(1/\vol R\bigr)} .
\]
\end{proposition}

\begin{proof}
Let $f\colon R\to[0,K]$ be a potential. Let $k\ge1$, $n=kd$ and $0<\delta<1$. For a point
$x=(x_1,\dots,x_k)$ of $R^k\subseteq[0,1)^n$, with $x_i\in R$, put $\Sigma(x)=f(x_1)+\dots+f(x_k)$. Since
$0\le\Sigma\le kK$, one of the sets
\[
  S_j=\bigl\{x\in R^k:\ j\delta^2/2\le\Sigma(x)<(j+1)\delta^2/2\bigr\},\qquad 0\le j\le2kK/\delta^2,
\]
has measure at least $(\vol R)^k\,\delta^2/(2kK+1)$. If $x,y,z\in S_j$ and $x+z\equiv2y\pmod1$, then summing
the inequality of Definition~\ref{def:potential} over the $k$ blocks of coordinates gives
$|x-z|^2\le\Sigma(x)+\Sigma(z)-2\Sigma(y)<\delta^2$, so all coordinates of $x-z$ have absolute value less
than $\delta$. This is the analogue of~\cite[Proposition~2.1]{EHPS}.

We now follow Green and Wolf~\cite{GreenWolf}, as at the end of~\cite[Section~3]{EHPS}. Let
$\delta=N^{-1/n}/4$. For a uniformly random $b\in[0,1)^n$ and a fixed $t\in\{1,\dots,N\}$, the point $tb$ is
uniformly distributed modulo~$1$, so it is congruent to a vector of maximum norm at most $\delta$ with
probability $(2\delta)^n$. Since $N(2\delta)^n=2^{-n}<1$, we can fix $b$ such that no $tb$ with
$1\le t\le N$ is congruent to such a vector. For $a\in[0,1)^n$ let $A_a$ be the set of
$t\in\{1,\dots,N\}$ with $a+tb\in S_j$ modulo~$1$. If $t_1+t_3=2t_2$ with $t_1>t_3$ in $A_a$, then the three
points $a+t_ib$ of $S_j$ satisfy $x+z\equiv2y$, so $(t_1-t_3)b$ is congruent to a vector of maximum norm
less than $\delta$, a contradiction. So $A_a$ is progression-free. For uniformly random $a$ the expected size
of $A_a$ is $N\vol S_j$, and hence
\[
  r_3(N)\;\ge\;\frac{N\,(\vol R)^{n/d}\,\delta^2}{2kK+1}\;=\;\frac{N\,(\vol R)^{n/d}\,N^{-2/n}}{16\,(2kK+1)} .
\]
Write $L=\log_2(1/\vol R)$, and take for $n$ a multiple of $d$ that differs from $\sqrt{2d\log_2N/L}$ by at
most $d$. Then $(n/d)L+(2/n)\log_2N=C_R\sqrt{\log_2N}+O(1)$, and $\log_2(16(2kK+1))=O(\log\log N)$,
which proves the proposition.
\end{proof}

The bound $K$ on the potential enters only the factor $1/(16(2kK+1))$, which disappears into the $o(1)$.
Hence a family of sets $R_\delta\subseteq R$ with $\vol R_\delta\to\vol R$ as $\delta\to0$, each of which
carries a potential, gives the bound with $C_R$ itself, however fast the bounds on the potentials grow.

\subsection*{The first floor}

\begin{lemma}\label{lem:Epot}
\begin{enumerate}
\item[(a)] The function $2\hf{u}^2$ is a potential on $F_1=\{u\in[0,1):\hf{u}<\frac14\}$.
\item[(b)] Let $\frac14\le\gamma\le\frac12$ and $0<\varepsilon\le\frac14$, and let $E(\gamma)_\varepsilon$ be
the set $E(\gamma)$ without the rows $\frac14-\varepsilon\le s<\frac14$ of its doubled part. Then
\[
  f_E(u,s)=q\,(\tfrac34-s)^2+2\hf{u}^2,\qquad q=\max\bigl(8,\,2+\tfrac1{2\varepsilon^2}\bigr),
\]
is a potential on $E(\gamma)_\varepsilon$, and $f_E\le\frac9{16}q+\frac12$.
\end{enumerate}
\end{lemma}

\begin{proof}
(a) If $x+z\equiv2y$ with $x,y,z\in F_1$, then by Lemma~\ref{lem:interval} with $I=[0,\frac14)$ the points
$x$ and $z$ lie in the same half of $[0,1)$ and $\hf{y}=\frac12(\hf{x}+\hf{z})$. So $x-z=\hf{x}-\hf{z}$,
and~\eqref{eq:midquad} gives the claim with equality.

(b) Let $x,y,z\in E(\gamma)_\varepsilon$ with $x+z\equiv2y\pmod1$. The levels satisfy
$s_x+s_z-2s_y\in\{-1,0,1\}$, since all levels lie in $[0,\frac34)$. Note that
$|x-z|^2=(u_x-u_z)^2+(s_x-s_z)^2<2$, and that the second term of $f_E$ changes by at least
$-4\hf{u_y}^2>-1$.

If $s_x+s_z=2s_y+1$, then property~(B) of $E(\gamma)$ (Lemma~\ref{lem:E}) is violated, so this case does
not occur.

If $s_x+s_z=2s_y-1$, put $a=\frac34-s_x$, $c=\frac34-s_z$ and $e=\frac34-s_y>0$. Then $a+c=2e+1$, so
$a^2+c^2-2e^2\ge\frac12(a+c)^2-2e^2=2e+\frac12>\frac12$, and the first term of $f_E$ contributes more than
$\frac q2\ge4$. This exceeds $1+|x-z|^2$.

If $s_x+s_z=2s_y$, the first term contributes $(q/2)(s_x-s_z)^2$ by~\eqref{eq:midquad}. Put
$D=\hf{u_x}+\hf{u_z}-2\hf{u_y}$. Since all values of $\hf{u}$ lie in $[0,\frac12)$ and $2u_y\equiv2\hf{u_y}$,
the number $D$ is $0$ or $\pm\frac12$, and as in the proof of Lemma~\ref{lem:interval}, $u_x$ and $u_z$ lie in
the same half of $[0,1)$ if and only if $D=0$. If $D=0$, then $u_x-u_z=\hf{u_x}-\hf{u_z}$, the second term
contributes exactly $(u_x-u_z)^2$, and $q/2\ge1$ finishes the proof. If $D\ne0$, we may assume that
$u_x\ge\frac12$. Then $x$ lies in the upper copy of the doubled part, and we write $p=\hf{u_x}<\frac14$ and
$r=u_z<\frac12$, so that $u_x-u_z=\frac12+p-r$.
\begin{itemize}
\item If $D=\frac12$, then $\hf{u_y}=\frac12(p+r-\frac12)\ge0$, so $r>\frac14$. The second term contributes
$2p^2+2r^2-(p+r-\frac12)^2=(p-r)^2+(p+r)-\frac14$, which exceeds
$(u_x-u_z)^2=(p-r)^2+(p-r)+\frac14$ by $2r-\frac12>0$.
\item If $D=-\frac12$, then $\hf{u_y}=\frac12(p+r+\frac12)\ge\frac14$, so $y$ lies in the single part and
$s_y\ge\frac14$, while $s_x<\frac14-\varepsilon$ because $x$ lies in the doubled part of
$E(\gamma)_\varepsilon$. Hence $s_z-s_x=2(s_y-s_x)>2\varepsilon$. The second term contributes
$(p-r)^2-(p+r)-\frac14$, which falls short of $(u_x-u_z)^2$ by $\frac12+2p<1$. The first term contributes
$(q/2)(s_x-s_z)^2\ge(s_x-s_z)^2+4(q/2-1)\varepsilon^2\ge(s_x-s_z)^2+1$, by the choice of $q$.
\end{itemize}
In all cases $f_E(x)+f_E(z)\ge2f_E(y)+|x-z|^2$. The bound on $f_E$ is clear.
\end{proof}

The rows that are removed in~(b) are needed for this potential: if $x$ lies in the upper copy of the doubled part just below
$s=\frac14$ and $z$ in the single part just above it, the mid-point $y$ can lie a quarter of a period further
along in $u$, and the second term of $f_E$ then falls short by almost $1$. The set of~\cite{EHPS} has a
similar gap.

\subsection*{The induction step}

\begin{theorem}\label{thm:steppot}
Let $B\subseteq[0,1)^d$ be a block, $F\subseteq[0,1)^{d-1}$ a measurable set, and $0\le\beta\le\frac14$.
Suppose that $B$ and $F$ carry potentials $f_B\le M_B$ and $f_F\le M_F$, and put $M=\max(M_B,M_F)$. For
$0<g,\eta\le\frac18$ let $B^+_{g,\eta}$ be the set $B^+(F,\beta)$ of Definition~\ref{def:step} without the
points of the high zone with $w<\frac14+g$ and without the points of the low zone with
$\frac14-\eta<\hf{t}<\frac14+\eta$. Then $B^+_{g,\eta}$ is a block, and the function $f$ given by
\[
  f(x',t,w)=f_B(x',t)+2w^2\quad\text{on the high zone},
\]
\[
  f(x',t,w)=C\,e^{-kw}+A\hf{t}^2+r(x',\hf{t})+2w^2\quad\text{on the low zone},
\]
where $r(x',\tau)=f_B(x',\tau)$ for $\tau<\frac14$, $r(x',\tau)=f_F(x')$ for $\tau\ge\frac14$, and
\[
  A=2+\frac{2M+d}{2\eta^2},\qquad K=2M_B+\frac A2+2M+d+4,\qquad k=\frac{2\ln3}g,\qquad C=3^{1/(2g)}K,
\]
is a potential on $B^+_{g,\eta}$ with $f\le C+A+M$.
\end{theorem}

\begin{proof}
Since $F$ carries a potential, it is reducible, so $B^+(F,\beta)$ is a block by Theorem~\ref{thm:step}, and
so is its subset $B^+_{g,\eta}$. The function $f$ is well defined, because $(x',t)\in B$ for a point of the
high zone, and $(x',\tau)\in B$ for a point of the low zone with $\tau=\hf{t}<\frac14$. We write $H$ and $L$
for the high and the low zone of $B^+_{g,\eta}$, so $w\in[\frac14+g,\frac34)$ on $H$ and $w\in[0,\frac14)$ on
$L$. On $H$ we have $f\le M_B+\frac98$. Let $x,y,z\in B^+_{g,\eta}$ with $x+z\equiv2y\pmod1$; note that
$|x-z|^2<d+1$. We distinguish five cases.

\emph{Case 1: $x,y,z\in H$.} Then $w_x+w_z-2w_y$ is an integer in $(-1,1)$, so $w_x+w_z=2w_y$. The points
$(x',t_x)$, $(y',t_y)$ and $(z',t_z)$ of $B$ satisfy the congruence, so $f_B$ contributes at least
$|x'-z'|^2+(t_x-t_z)^2$, and by~\eqref{eq:midquad} the term $2w^2$ contributes $(w_x-w_z)^2$.

\emph{Case 2: $y\in H$, and $x\in L$ or $z\in L$.} Say $x\in L$. Then $f(x)\ge Ce^{-kw_x}>Ce^{-k/4}=K$, while
$2f(y)+|x-z|^2<2M_B+\frac94+d+1<K$.

\emph{Case 3: $x,y,z\in L$.} Then $w_x+w_z=2w_y$. Write $\tau=\hf{t}$. The points $(\tau,w)$ of the low zone
lie in the region $\Omega_\beta$ of Lemma~\ref{lem:chart}, with $(\tau,w)$ in the roles of $(\sigma,\tau)$. By
that lemma, applied to the coordinates $t$ and $w$, we have $\tau_x+\tau_z=2\tau_y$, and $t_x$ and $t_z$ lie in
the same half of $[0,1)$, so $t_x-t_z=\tau_x-\tau_z$. The term $Ce^{-kw}$ is convex in $w$ and contributes at least $0$, the term $2w^2$
contributes $(w_x-w_z)^2$, and $A\tau^2$ contributes $(A/2)(\tau_x-\tau_z)^2$. For the term $r$ there are
three possibilities.
\begin{itemize}
\item If $\tau_x,\tau_z<\frac14$, then $\tau_y<\frac14$, and $(x',\tau_x)$, $(y',\tau_y)$, $(z',\tau_z)$ are
points of $B$ that satisfy the congruence. So $r$ contributes at least $|x'-z'|^2+(\tau_x-\tau_z)^2$.
\item If $\tau_x,\tau_z\ge\frac14$, then $\tau_y\ge\frac14$, and $x',y',z'\in F$ satisfy $x'+z'\equiv2y'$.
So $r$ contributes at least $|x'-z'|^2$, and $(A/2)(\tau_x-\tau_z)^2\ge(\tau_x-\tau_z)^2$.
\item Otherwise one of $\tau_x,\tau_z$ is at most $\frac14-\eta$ and the other at least $\frac14+\eta$, so
$|\tau_x-\tau_z|\ge2\eta$. The term $r$ contributes at least $-2M$, and
$(A/2-1)(\tau_x-\tau_z)^2\ge2(A-2)\eta^2=2M+d$, which exceeds $2M+|x'-z'|^2$.
\end{itemize}
In each case $f(x)+f(z)-2f(y)\ge|x'-z'|^2+(t_x-t_z)^2+(w_x-w_z)^2=|x-z|^2$.

\emph{Case 4: $y\in L$, and one of $x,z$ lies in $L$ and the other in $H$.} Say $x\in L$ and $z\in H$. Then
$w_x+w_z-2w_y$ is an integer in $(-\frac14,1)$, so $w_z=2w_y-w_x$. Since $w_z\ge\frac14+g$ and
$w_x<\frac14$, we get $w_y-w_x>\frac g2$, so $Ce^{-kw_x}>3Ce^{-kw_y}$, and the exponential terms contribute
more than $Ce^{-kw_y}>Ce^{-k/4}=K$. The other terms of $f(x)$ and $f(z)$ are non-negative, the other terms
of $2f(y)$ are at most $2(\frac A4+M+\frac18)$, and $|x-z|^2<d+1$. The choice of $K$ finishes this case.

\emph{Case 5: $y\in L$ and $x,z\in H$.} Then $w_x+w_z=2w_y+1$, and property~(B) of $B^+_{g,\eta}$ shows that
$x+z\not\equiv2y$. So this case does not occur.

Finally, $f\le C+\frac A4+M+\frac18\le C+A+M$ on $L$, and $f\le M_B+\frac98<C$ on $H$.
\end{proof}

The two gaps are needed for this proof. Without the gap above $w=\frac14$ the exponential term cannot
separate the low zone from the high zone in Case~4. Without the gap at $\hf{t}=\frac14$ the potentials of
the inherited and the free rows, which are unrelated, cannot be combined in Case~3. In special cases the
second gap is not needed (Theorem~\ref{thm:Mpot}), and we do not claim that the gaps are needed for every
potential. The volume of $B^+(F,\beta)$ is a linear function of the
slice measure of $B$ and of $\vol F$, by~\eqref{eq:volplus}, with weights $\omega_h+2\omega_\ell\le\frac34$
and $\beta^2$. So if $B'\subseteq B$ is a measurable block and $F'\subseteq F$ is measurable, then
\begin{equation}\label{eq:lipschitz}
  \vol\bigl(B^+(F,\beta)\setminus B'^+_{g,\eta}(F',\beta)\bigr)\;\le\;\tfrac34\vol(B\setminus B')
  +\beta^2\vol(F\setminus F')+g+\eta ,
\end{equation}
where $g$ bounds the volume of the removed levels of the high zone, and $\eta$ that of the removed rows of
the low zone, which have measure at most $4\eta$ at every level $w<\frac14$.

\begin{corollary}\label{cor:towerpot}
Let $T_2,T_3,\dots$ be the tower of Definition~\ref{def:tower}, with any parameters. For every $d\ge2$ and
every $\delta>0$ there is a block $T_d'\subseteq T_d$ with $\vol(T_d\setminus T_d')<\delta$ that carries a
potential.
\end{corollary}

\begin{proof}
For $d=2$ take $T_2'=E(\gamma)_\varepsilon$ with $\varepsilon<2\delta$ (Lemma~\ref{lem:Epot}(b)); the removed
rows have measure at most $\frac12$ each. The free sets of the tower are $F_1$, which carries a potential by
Lemma~\ref{lem:Epot}(a), $E(\frac5{12})$, which may be replaced by $E(\frac5{12})_\varepsilon$, and the
floors $T_{d-1}$, which may be replaced by $T_{d-1}'$. Applying Theorem~\ref{thm:steppot} floor by floor and
using~\eqref{eq:lipschitz}, with all thinnings small enough, gives the claim.
\end{proof}

\subsection*{\texorpdfstring{The sets $R_\beta$ and $T\cup M^*$}{The sets R\_beta and T u M*}}

The sets $R_\beta$ and $T\cup M^*$ of Section~\ref{sec:dim3} are not blocks, so Theorem~\ref{thm:steppot} does not
apply to them, and their removal orders are not level by level. We give explicit potentials instead, for both
dimensions and all $0<\beta\le\frac14$ at once. The low zone and the notch share one formula in the chart
coordinates of Section~\ref{sec:dim3}, with a term that follows the lines of the fan.

\begin{theorem}\label{thm:Mpot}
Let $d\in\{2,3\}$, $0<\beta\le\frac14$ and $\frac14\le\gamma\le\frac12$, and let $0<g,\varepsilon\le\frac18$,
$0<\varepsilon_0<\min(\beta,\frac18)$ and $0<\lambda<1$. With the notation of Definition~\ref{def:Mstar} let
\[
  N'=\bigl\{(s,t)\in N:\ s(1+8\tau)<\lambda(\beta-\tau)(1-8\tau)\ \text{ and }\
  \tau\le\min(\beta,\tfrac18)-\varepsilon_0\bigr\},\qquad\tau=t-\tfrac12,
\]
and $M'=\{(u,s,t)\in M^*:\ (s,t)\in N'\}$. For $d=2$ let $P'$ be the union of $L$, $N'$ and the set $H'$ of the
points of $H$ with $t\ge\frac14+g$, and put $f_B(s)=2s^2$ and $M_B=\frac12$. For $d=3$ let $P'$ be the union of
$M'$ and the set $T'$ obtained from $T=T(\beta,\gamma)$ by removing the points of the high zone with
$t<\frac14+g$ or with $\frac14-\varepsilon\le s<\frac14$, let $H'$ be the high zone of $T'$, and let $f_B=f_E$ and
$M_B=\frac9{16}q+\frac12$ with $q=\max(8,\,2+1/(2\varepsilon^2))$, as in Lemma~\ref{lem:Epot}(b). On the low
zone and on the notch use the chart coordinates $\sigma=\hf{s}$ and $\tau=\hf{t}$, put $v=\hf{u}$ for $d=3$ and
$v=0$ for $d=2$, put $\nu=1-8\tau$, and let $\Psi(\sigma,\tau)=A\,\nu\,e^{-k_2\sigma/\nu}$ for $\tau<\frac18$ and
$\Psi(\sigma,\tau)=0$ for $\tau\ge\frac18$. Then the function
\[
  f=f_B+2t^2\quad\text{on }H',
\]
\[
  f=E_0+C\,e^{-k_1\tau}+\Psi(\sigma,\tau)+2(\sigma^2+\tau^2+v^2)\quad\text{on the low zone and on the notch of }P',
\]
with
\begin{gather*}
  E_0=2M_B+\tfrac94+d,\qquad K_1=E_0+d+\tfrac32,\qquad k_1=\frac{2\ln3}g,\qquad C=3^{1/(2g)}K_1,\\
  K_2=K_1+2C,\qquad \nu_0=8\varepsilon_0,\qquad k_2=\frac{2\ln(4/\nu_0)}{(1-\lambda)\,\varepsilon_0},\qquad
  A=\frac{2K_2}{\nu_0}\,e^{\lambda\beta k_2},
\end{gather*}
is a potential on $P'$. As $g,\varepsilon,\varepsilon_0\to0$ and $\lambda\to1$, the volume of $P'$ tends to
$\vol R_\beta$ for $d=2$ and to $\vol(T\cup M^*)$ for $d=3$.
\end{theorem}

\begin{proof}
For $d=2$ the function $f_B$ is the potential $2s^2$ of the interval $[0,\frac12)$, and for $d=3$ it is the
potential $f_E$ of Lemma~\ref{lem:Epot}(b) on $E(\gamma)_\varepsilon$. On $H'$ we have $t<\frac34$, so
$f<M_B+\frac98$. We write $L'$ for the union of the low zone and the notch of $P'$. On $L'$ we have
$\sigma<\frac12$, $\tau<\frac14$ and $v<\frac14$, so $4(\sigma^2+\tau^2+v^2)<\frac32$, and on the notch of $P'$ we
have $\nu\ge\nu_0$. We also need that $\Psi$ is convex. For $\tau<\frac18$ it is the perspective
$(\sigma,\nu)\mapsto\nu G(\sigma/\nu)$ of the convex function $G(\theta)=Ae^{-k_2\theta}$, which is convex for
$\nu>0$, increasing in $\nu$, and tends to $0$ as $\nu\to0$. So $\Psi$, as a function of $\sigma$ and
$\max(\nu,0)$, is convex. Let $x,y,z\in P'$ with $x+z\equiv2y\pmod1$; note that $|x-z|^2<d$.

\emph{Case 1: $x,y,z\in H'$.} This is Case~1 of the proof of Theorem~\ref{thm:steppot}, with the potential
$f_B$.

\emph{Case 2: $y\in H'$, and $x$ or $z$ lies in $L'$.} Then $f(x)+f(z)\ge E_0=2(M_B+\frac98)+d$, which exceeds
$2f(y)+|x-z|^2$.

\emph{Case 3: $y\in L'$ and $x,z\in H'$.} This is impossible, by Case~1 of the proof of
Theorem~\ref{thm:Mstar}, which does not use the order.

\emph{Case 4: $x,y,z\in L'$.} By Lemma~\ref{lem:chart}, the chart point of $y$ is the mid-point of the chart
points of $x$ and $z$, and $s_x-s_z=\sigma_x-\sigma_z$ and $t_x-t_z=\tau_x-\tau_z$. For $d=3$,
Lemma~\ref{lem:interval} with $I=[0,\frac14)$ shows that $u_x$ and $u_z$ lie in the same half of $[0,1)$ and
that $v_y=\frac12(v_x+v_z)$. So $|x-z|^2=(v_x-v_z)^2+(\sigma_x-\sigma_z)^2+(\tau_x-\tau_z)^2$, which is exactly the
contribution of $2(\sigma^2+\tau^2+v^2)$, by~\eqref{eq:midquad}. The terms $Ce^{-k_1\tau}$ and $\Psi$ are convex and
contribute at least $0$.

\emph{Case 5: $y\in L'$, and one of $x,z$ lies in the low zone and the other in $H'$.} Say $x$ lies in the low
zone and $z\in H'$. Then $t_z\equiv2\tau_y-\tau_x\pmod1$ with $2\tau_y-\tau_x\in(-\frac14,\frac12)$, so
$t_z=2\tau_y-\tau_x\ge\frac14+g$. Hence $\tau_y\ge\frac18+\frac g2$, so $\Psi(y)=0$, and $\tau_y-\tau_x>\frac g2$, so
$Ce^{-k_1\tau_x}-2Ce^{-k_1\tau_y}>Ce^{-k_1\tau_y}\ge Ce^{-k_1/4}=K_1$. Therefore
\[
  f(x)+f(z)-2f(y)-|x-z|^2\;>\;K_1-E_0-\tfrac32-d\;=\;0 .
\]

\emph{Case 6: $y\in L'$, and one of $x,z$ lies in the notch and the other in $H'$.} Say $x$ lies in the notch
and $z\in H'$. Then $t_x=\tau_x+\frac12$ and $2t_y\equiv2\tau_y$, so $t_z\equiv2\tau_y-\tau_x+\frac12\pmod1$ with
$2\tau_y-\tau_x\in(-\frac18,\frac12)$, and since $t_z\in[\frac14,\frac34)$ we get $t_z=2\tau_y-\tau_x+\frac12$.
Similarly $s_z\equiv2\sigma_y-\sigma_x$ with $2\sigma_y-\sigma_x\in(-\frac14,1)$, where a negative value
would put $s_z$ above $\frac34$, so $s_z=2\sigma_y-\sigma_x$. The cut $s_z\ge t_z-\frac12+\beta$ of the high zone
gives $\sigma_y-\tau_y\ge w^*$, where $w^*=\frac12(\sigma_x-\tau_x+\beta)$. Hence
$\Psi(y)\le Ae^{-k_2\sigma_y}\le Ae^{-k_2w^*}$, also when $\tau_y\ge\frac18$. With $\nu=\nu_x$ and
$\theta=\sigma_x/\nu_x$ we have $\tau_x=\frac{1-\nu}8$ and
\[
  w^*-\theta\;=\;\tfrac12\bigl((\beta-\tau_x)-\theta(2-\nu)\bigr).
\]
The first inequality that defines $N'$ reads $\theta(2-\nu)<\lambda(\beta-\tau_x)$. So
$w^*-\theta>\frac12(1-\lambda)(\beta-\tau_x)\ge\frac12(1-\lambda)\varepsilon_0$, and $\theta<\lambda\beta$ because
$2-\nu=1+8\tau_x\ge1$. Since $\nu\ge\nu_0$ and $e^{-k_2(1-\lambda)\varepsilon_0/2}=\nu_0/4$, this gives
$2\Psi(y)\le2Ae^{-k_2\theta}e^{-k_2(w^*-\theta)}<(\nu_0/2)\,Ae^{-k_2\theta}$, while
$\Psi(x)=A\nu e^{-k_2\theta}\ge\nu_0Ae^{-k_2\theta}$. Hence
$\Psi(x)-2\Psi(y)>(\nu_0/2)\,Ae^{-k_2\lambda\beta}=K_2$, and
\[
  f(x)+f(z)-2f(y)-|x-z|^2\;>\;K_2-E_0-2C-\tfrac32-d\;=\;0 .
\]

So $f$ is a potential. The removed parts of $L\cup H$ and of $T$ have volume at most $g+\varepsilon$, and the parts of
the notch that are not in $N'$ or $M'$ shrink to null sets as $\varepsilon_0\to0$ and $\lambda\to1$.
\end{proof}

Only Case~6 needs $\lambda<1$: on the curved boundary of the notch the difference $w^*-\theta$ vanishes, which
is the tightness noted in Remark~\ref{rem:sharp}. Since the free set in Theorem~\ref{thm:steppot} needs a
potential but not property~(B), the thinned sets of the theorem can serve as free sets, for instance the set
$P'$ for $d=2$ and $\beta=\frac1{10}$ in a fourth floor.

\begin{proof}[Proof of Theorem~\ref{thm:behrend} and of the constants $C$ in Theorem~\ref{thm:dim3}]
Let $R$ be a set in Table~\ref{tab:constants}, of dimension $d$. By the remark after
Proposition~\ref{prop:behrend} it suffices to find, for every $\delta>0$, a subset of $R$ of volume at least
$\vol R-\delta$ that carries a potential. For $[0,\frac12)$ this is the interval itself, for $E(\frac5{12})$ it
is Lemma~\ref{lem:Epot}(b), for $R_{1/10}$, $T_3$ and $T(\frac7{50},\frac25)\cup M^*$ it is
Theorem~\ref{thm:Mpot}, and for the floors $T_d$, $d\ge4$, it is Corollary~\ref{cor:towerpot}. For
$\widetilde T_4$ we apply Theorem~\ref{thm:steppot} to a thinned block as in Corollary~\ref{cor:towerpot} and to
the thinned set $P'\subseteq R_{1/10}$ of Theorem~\ref{thm:Mpot} as the free set, and use~\eqref{eq:lipschitz}.
Proposition~\ref{prop:behrend} then gives the bound with $C=2\sqrt{(2/d)\log_2(1/\vol R)}$, and the volumes of
Theorem~\ref{thm:dim3} and Proposition~\ref{prop:higher} give the values in the table.

For the first statement let $C_\infty=2\sqrt{\log_2(r_\infty^{-2})}$ and $\varepsilon>0$. By
Theorem~\ref{thm:main}(a) there is a floor $T_d$ with $2\sqrt{(2/d)\log_2(1/\vol T_d)}<C_\infty+\varepsilon$, so
$r_3(N)\ge N\cdot2^{-(C_\infty+\varepsilon+o(1))\sqrt{\log_2N}}$. Since this holds for every $\varepsilon>0$, it
holds with $C_\infty+o(1)$ in the exponent. Finally $r_\infty\ge0.55499308$ gives
$C_\infty\le2\sqrt{\log_2(0.55499308^{-2})}<2.60685$.
\end{proof}

\begin{remark}\label{rem:sharpBehrend}\label{rem:Ksize}
The bounds on our potentials are enormous, and they grow with every floor and as the thinning shrinks. This
does not matter for Theorem~\ref{thm:behrend}, because the bound enters only the factor $1/(16(2kK+1))$ in the
proof of Proposition~\ref{prop:behrend}; it is the reason why we fix $d$ and the thinning while $N\to\infty$,
whereas in~\cite{EHPS} the gap shrinks with the dimension. We do not know whether $R_\beta$ or $T\cup M^*$
itself, without thinning, carries a potential. Numerical experiments suggest that they do not: near the curved
boundary of the notch the inequality forces a potential to halve again and again along level~$0$. This is a
heuristic observation, not a proof.
\end{remark}

\section{Upper bounds}\label{sec:upper}

We collect upper bounds for reducible sets and for the sets that the induction step can produce. We write $M(m,d)$ for the largest size of a reducible subset of $(\Zm)^d$,
\[
  v_d=\sup\{\vol S:\ S\subseteq[0,1)^d\text{ measurable and reducible}\},
\]
and $\kappa^*=\sup_{d\ge1}v_d^{1/d}$. Products of reducible sets are reducible (for a finite subset of a product, add the functions $h$ of
Remark~\ref{rem:eps24}, whose construction works in every abelian group, for its two projections; a point where
the sum is largest is a non-mid-point). So $v_{d+e}\ge v_dv_e$, and
therefore $\kappa^*=\lim_{d\to\infty}v_d^{1/d}$ by Fekete's lemma. Since a set that carries a potential is
reducible (Section~\ref{sec:behrend}), all bounds of this section also hold for such sets.

\begin{proposition}[Reflection bound]\label{prop:reflection}
Let $S\subseteq(\Zm)^d$ be reducible. Then $|S|\le\frac12(m^d+1)$ if $m$ is odd, and
$|S|\le\frac12(m^d+2^d)$ if $m$ is even.
\end{proposition}

\begin{proof}
We may assume that $S$ is not empty. Let $y\in S$ be a point that is not a mid-point of $S$. If $x$ and $2y-x$
both lie in $S$, then $x+(2y-x)=2y$, hence $x=2y-x$ by the choice of $y$, that is, $2(x-y)=0$. This has one
solution $x$ if $m$ is odd and $2^d$ solutions if $m$ is even. So $S$ and its reflection $2y-S$ have at most
$1$, respectively $2^d$, common points, and $2|S|-|S\cap(2y-S)|\le m^d$.
\end{proof}

For $d=1$ and odd $m$ the bound is attained by the interval $\{0,1,\dots,\frac{m-1}{2}\}$. Applied to the points
of $S$ on a line $\{a+kv:k\in\Zm\}$ whose direction $v$ has order $m$, the proposition shows that a reducible
set has at most $\frac12(m+1)$ points on such a line if $m$ is odd, and at most $\frac12(m+2)$ if $m$ is even.

\begin{proposition}[Small moduli]\label{prop:small}
\begin{enumerate}
\item[(a)] A subset of $\mathbb{F}_3^d$ is reducible if and only if it is a cap. Hence $M(3,d)$ is the largest
size of a cap in $\mathrm{AG}(d,3)$: $M(3,d)=2,4,9,20,45,112$ for $d=1,\dots,6$~\cite{Pellegrino,EFLS,Potechin}.
\item[(b)] $M(5,1)=3$ and $M(5,2)=10$.
\end{enumerate}
\end{proposition}

\begin{proof}
(a) A reducible set has at most two points on every line, so it is a cap. Conversely, in $\mathbb{F}_3^d$ three
distinct points $x,y,z$ satisfy $x+z=2y$ if and only if they are collinear, so a cap has no mid-points at all.
(b) Since subsets of reducible sets are reducible, all reducible subsets of $\mathbb{F}_5^2$ can be enumerated by
adding one point at a time. This enumeration finds $88000$ reducible sets with $10$ points and none with $11$.
\end{proof}

\begin{proposition}[Polynomial method]\label{prop:EG}
Let $p$ be an odd prime and let $S\subseteq\mathbb{F}_p^d$ be reducible. Then $|S|^{1/d}\le J_p\,p$, where
\[
  J_p=\frac1p\,\min_{0<t<1}\frac{1-t^p}{(1-t)\,t^{(p-1)/3}}
\]
is the constant of Ellenberg and Gijswijt~\cite{EG}. We have $J_3<0.9184$, $J_5<0.8924$, $J_7<0.8795$, and $J_p$
decreases to $J_\infty=0.84143\ldots<0.8415$ as $p\to\infty$~\cite{BCCGNSU}.
\end{proposition}

\begin{proof}
By Proposition~\ref{prop:eps}, $\mathbb{F}_p^n$ contains a progression-free set of size $|S|^{n/d-O(\log n)}$. We
compare this with $r_3(\mathbb{F}_p^n)\le(J_pp)^n$~\cite{EG} and let $n\to\infty$.
\end{proof}

\begin{remark}\label{rem:finite}
For a fixed odd prime $p$ the reflection bound reads $|S|^{1/d}/p\le(\frac12(1+p^{-d}))^{1/d}$. It is better than
Proposition~\ref{prop:EG} in low dimensions and worse in high dimensions, and the dimension where this changes
depends on $p$; in the limit $p\to\infty$ the reflection bound is the stronger one exactly for $d\le4$. For a
fixed modulus $m$ the limit $\rho(m)=\lim_{n\to\infty}r_3((\Zm)^n)^{1/n}$ exists by Fekete's lemma, since $r_3$ is
supermultiplicative. Progression-free sets are reducible, and the proof above gives $M(m,d)^{1/d}\le\rho(m)$, so
$\lim_{d\to\infty}M(m,d)^{1/d}=\rho(m)$. Hence the optimal asymptotic exponential bounds for progression-free
sets and for reducible sets coincide. Hametner and Tyrrell~\cite{HT} proved that for odd primes $p$ and
progression-free sets $A\subseteq\mathbb{F}_p^d$ the inequality $|A|^{1/d}<\rho(p)$ is strict in every dimension; we do not know whether
it is strict for reducible sets.
\end{remark}

\begin{proposition}\label{prop:torus}
Let $S\subseteq[0,1)^d$ be measurable and reducible. Then
\begin{enumerate}
\item[(a)] $\vol S\le\frac12$;
\item[(b)] $\vol S\le M(m,d)/m^d$ for every $m\ge2$, and $\vol S\le J_\infty^{\,d}$;
\item[(c)] $\vol S\le v_{d-k}\cdot\lambda^*_k(\pi(S))$ for $1\le k<d$, where $\pi$ is the projection to any $k$ of
the coordinates and $\lambda^*_k$ is the outer Lebesgue measure. In particular $v_d\le v_{d-1}$.
\end{enumerate}
\end{proposition}

\begin{proof}
By Lemmas~\ref{lem:transfer} and~\ref{lem:average}, for every $m$ there is a reducible subset of $(\Zm)^d$ with
at least $\vol(S)\,m^d$ elements. This gives the first part of~(b). Part (a) follows from
Proposition~\ref{prop:reflection} with odd $m\to\infty$, and the second part of~(b) from
Proposition~\ref{prop:EG} with primes $p\to\infty$. For (c), fix the values of $k$ coordinates. If three points
$x,y,z$ of the corresponding section satisfy $x+z\equiv2y$ in the remaining $d-k$ coordinates, then they
satisfy it in all coordinates. So every section is a reducible subset of $[0,1)^{d-k}$, almost every section
is measurable, and (c) follows from Fubini's theorem.
\end{proof}

\begin{corollary}\label{cor:vd}
We have $v_1=\frac12$, the sequence $v_d$ is non-increasing, and
\[
  v_d\le\min\Bigl(v_{d-1},\ \inf_{m\ge2}\frac{M(m,d)}{m^d},\ J_\infty^{\,d}\Bigr)\quad(d\ge2),\qquad
  \kappa^*\le J_\infty .
\]
In particular $v_2\le\frac{10}{25}=\frac25$, $v_3\le\frac9{27}=\frac13$, $v_4\le\frac{20}{81}$,
$v_5\le\frac{45}{243}$ and $v_6\le\frac{112}{729}$, by Proposition~\ref{prop:small}.
\end{corollary}

In
dimension~$2$ the largest area of a reducible subset of the torus lies between $0.2943$
(Theorem~\ref{thm:dim3}) and $\frac25$, and we expect that it is far below $\frac25$. With the constructions of
this paper the state of knowledge is as follows.

\begin{table}[ht]
\centering
\small
\renewcommand{\arraystretch}{1.2}
\begin{tabular}{@{}lllll@{}}
\toprule
$d$ & $v_d\ge$ & $v_d\le$ & $v_d^{1/d}\ge$ & $v_d^{1/d}\le$\\
\midrule
$1$ & $1/2$ & $1/2$ & $0.5$ & $0.5$\\
$2$ & $0.2943$ & $2/5$ & $0.5425$ & $0.6325$\\
$3$ & $0.1658$ & $1/3$ & $0.5494$ & $0.6934$\\
$4$ & $0.0918$ & $20/81$ & $0.5505$ & $0.7050$\\
$5$ & $0.0508$ & $45/243$ & $0.5511$ & $0.7138$\\
$6$ & $0.0282$ & $112/729$ & $0.5517$ & $0.7319$\\
$d\ge7$ & $\vol T_d$ & $\min(112/729,\,J_\infty^{\,d})$ & $(\vol T_d)^{1/d}$ & $0.8415$\\
$d\to\infty$ & & & $0.554993$ & $0.8415$\\
\bottomrule
\end{tabular}
\caption{Bounds for the largest volume $v_d$ of a measurable reducible subset of $[0,1)^d$. The lower
bounds come from Theorem~\ref{thm:dim3} in dimensions $2$, $3$ and $4$ and from the tower in higher dimensions.
Lower bounds are rounded down and upper bounds are rounded up.}
\label{tab:vd}
\end{table}

\subsection*{Bounds for the induction step}

The sets produced by the induction step satisfy stronger bounds.

\begin{proposition}\label{prop:stepbound}
\begin{enumerate}
\item[(a)] For $d\ge2$, every measurable block $B\subseteq[0,1)^d$ satisfies $\vol B\le\frac34\,v_{d-1}$.
\item[(b)] Let $B$, $F$ and $\beta$ be as in Definition~\ref{def:step}, with Borel sets $B$ and $F$, and let
$\mu=\max\bigl(\sup_{\tau<1/4}\vol B_\tau,\ \vol F\bigr)$ be the largest volume of a row. Then
\[
  \vol B^+\le\tfrac12\vol B+\tfrac18\mu,\qquad \sup_{w<1/4}\vol\,(B^+)_w\le\tfrac12\mu .
\]
\item[(c)] Let $T_{d+1}=T_d^+(T_{d-1},\beta_{d+1})$ for $d\ge d_0$ be a tower of Borel sets with arbitrary
parameters $\beta_d$. Then
\[
  \limsup_{d\to\infty}(\vol T_d)^{1/d}\;\le\;\frac{1+\sqrt3}4\;=\;0.68301\ldots .
\]
\end{enumerate}
\end{proposition}

\begin{proof}
(a) The levels of a block are below $\frac34$, and its slices are reducible.

(b) The high zone occupies the levels $w\in[\frac14,\frac34)$, and its slice at each level is a subset of $B$.
The slice of the low zone at a level $w<\frac14$ consists of the rows with $\hf{t}$ in a window of length at
most $\frac14$, each taken twice, so its volume is at most $2\cdot\frac14\cdot\mu$. Integrating over
$w\in[0,\frac14)$ gives the first inequality.

(c) Write $a_d=\vol T_d$ and $s_d=\sup_{w<1/4}\vol\,(T_d)_w$. By (b),
\[
  a_{d+1}\le\tfrac12a_d+\tfrac18\max(s_d,a_{d-1}),\qquad s_{d+1}\le\tfrac12\max(s_d,a_{d-1}).
\]
Let $\lambda=(1+\sqrt3)/4$, so that $\lambda^2=\frac12\lambda+\frac18$ and $\lambda>\frac12$, and choose $D$ so
large that $a_d\le D\lambda^d$ and $s_d\le D\lambda^{d-1}$ for $d=d_0-1$ and $d=d_0$. If these bounds hold for
$d-1$ and $d$, then $a_{d+1}\le\frac12D\lambda^d+\frac18D\lambda^{d-1}=D\lambda^{d+1}$ and
$s_{d+1}\le\frac12D\lambda^{d-1}\le D\lambda^d$. So they hold for all $d$ by induction.
\end{proof}

Part~(c) concerns towers in which every floor has the floor two dimensions lower as its free set, as in
Definition~\ref{def:tower}; it does not cover free sets that are chosen independently on every floor, and it
does not cover sets that are not removed level by level, such as $R_\beta$ and $T\cup M^*$. For Behrend's
construction, $\vol R\le J_\infty^{\,d}$ shows that the method of Proposition~\ref{prop:behrend} cannot give a
constant $C$ below $2\sqrt{2\log_2(1/J_\infty)}>1.41$, and part~(c) limits the asymptotic growth of such towers,
which corresponds to $C=2\sqrt{2\log_2(4/(1+\sqrt3))}>2.09$. These are limits of the method; they say nothing
about the true order of $r_3(N)$.

\section{Concluding remarks}\label{sec:concluding}

\subsection*{How the constructions were found}
The three-dimensional set $T$ was found by a computer search in a restricted class, followed by proofs.
We considered sets fibred over the $(s,t)$-plane whose $u$-fibres are either two copies of an interval
($\hf{u}<\rho$) or a single interval $[c,\frac12)$. Once the pattern of supports is fixed, the sufficient
conditions that we used for reducibility with the bottom-up removal order are linear in the fibre lengths, so
the largest set of this restricted class with a given pattern is the solution of a linear programme, and one can search over patterns with the linear
programme inside. At resolution $\frac1{32}$ the optimum had a clean piecewise linear shape. Reading off
its slopes and breakpoints gave a first version of volume $35/216$, which had an additional cut that the
proof does not need; removing the cut and re-balancing the parameters gave $T(\frac18,\frac25)$.
Simulated annealing at $m=15$, seeded with discretisations of these sets, then polished their boundaries,
and comparing the polished sets with their seeds layer by layer showed where the formal sets had slack.
This led to the region $M^*$. Its section at $u=0$ is the notch $N$ of the two-dimensional set $R_\beta$, and the
proof of Theorem~\ref{thm:Mstar} works in both dimensions.

The potentials of Section~\ref{sec:behrend} were found in the same way. Linear programmes for potentials
on grids have optimal solutions with a clear shape: on the low zone they decrease
exponentially with the level, and on $M^*$ they follow the lines of the fan and halve along the chains
described in Remark~\ref{rem:sharpBehrend}. Writing these shapes as formulas gave Theorems~\ref{thm:steppot}
and~\ref{thm:Mpot}.

The induction step came from asking which features of the passage from dimension~$2$ to dimension~$3$
are general. There are three. First, a product uses the same half of every line in the new direction, so
extra volume needs more than half of the levels, with different slices at different levels. Second, the
only new obstructions are triples that wrap around once in the level coordinate. They are excluded when
the sums of two high slices miss the doubles of the low slice. The map $x\mapsto2x$ of the torus preserves
measure in the sense that the preimage of a set has the measure of the set, so the low slice, which may be any
reducible part of the preimage of that gap, can be as large as the gap. Third, the low rows of a block already avoid
the sums of its high rows, so they can be used again one dimension higher. The third point is what makes
the step repeatable, and it is property~(B).

\subsection*{Open problems}
\begin{enumerate}
\item \emph{The limit of the tower.} What is the best constant that the induction step can give, with
optimal cuts, windows, free sets and parameters? With straight cuts and windows and $\beta=\frac16$ on all
high floors the growth rate is $r_\infty=0.554993\ldots$ (Theorem~\ref{thm:main}). It is open how much can be
gained by other sequences of parameters, for instance periodic ones, by convex cuts with the windows of
Remark~\ref{rem:step}(b), or by free sets that are chosen for their effect on later floors and not for their
volume.
\item \emph{Refinements at higher floors.} In dimensions~$2$ and~$3$ the notches and the fan order raised the
constants beyond the tower. Is there a version of the induction step that allows removal orders that do not
proceed level by level, so that analogues of $M^*$ can be added on every floor?
Proposition~\ref{prop:Mnext} shows that the step as defined here cannot simply be applied to $T_3\cup M^*$.
\item \emph{Other mechanisms.} The tower gains one coordinate at a time. Can two level coordinates
interact? Can an inherited row be larger than the corresponding slice of the block? Only pairs of points
above the cut constrain it.
\item \emph{Dimensions two and three.} What is the largest area of a reducible subset of the two-dimensional
torus? It lies between $0.2943$ and $\frac25$. Are the notches the largest regions of their kind
(Remark~\ref{rem:sharp}), is there a larger region with a different pivot or a different order, and is
$T\cup M^*$ optimal among fibred sets? Our searches were not exhaustive.
\item \emph{Upper bounds.} The constructions give constants up to $0.554993$, while the best upper bounds
are $(1+\sqrt3)/4$ for the growth of the towers of Proposition~\ref{prop:stepbound}(c) and $J_\infty<0.8415$ for
reducible sets in general; the latter rests on the bound of Ellenberg and Gijswijt for odd primes. Can the
bound $(1+\sqrt3)/4$ be improved, or extended to arbitrary free sets? Is $\kappa^*<J_\infty$?
\item \emph{Potentials.} Can every measurable reducible subset of a torus be approximated in volume by sets
that carry potentials? For the sets of this paper this is the case, but the sharp sets $R_\beta$ and
$T\cup M^*$ may carry none (Remark~\ref{rem:sharpBehrend}). Can the potentials of Section~\ref{sec:behrend} be
chosen with bounds that grow slowly enough to let the thinning shrink with the dimension, as in~\cite{EHPS}?
This would improve the $o(1)$ term in Theorem~\ref{thm:behrend}.
\end{enumerate}

\subsection*{Computations}
Four results of this paper rest on computations, all in exact integer arithmetic. The programmes are supplied
with the paper as ancillary files, with the commands to compile and run them: \file{tower\_exact.c} for the
volumes in Propositions~\ref{prop:T4improved} and~\ref{prop:higher}, \file{spectral\_cert.c} for
Proposition~\ref{prop:cert}, \file{core\_check.c} for Proposition~\ref{prop:Mnext}, and \file{m5\_search.c} for
Proposition~\ref{prop:small}(b). In addition we tested many discretisations of all sets of this paper in
$(\Zm)^d$ by removing non-mid-points, and the potentials on grids. Searches by simulated annealing, seeded
with discretisations of our sets, and linear programmes for the sizes of the fibres and for potentials led
to the constructions, as described above; they found no further region of positive volume, but they are
heuristic and prove nothing.

\appendix
\section{A more general removal rule}\label{app:rule}

One may ask whether a more permissive removal rule could give larger sets. Instead of removing only points
that are not the middle term of any progression, one could remove every point that is missing from at
least one of the three positions of the non-constant progressions in the current set. For three-term
progressions and odd $m$ this reduces exactly the same sets.

A \emph{non-constant $3$-progression} in a set $S$ is an ordered triple $(x,y,z)$ of points of $S$, not all
equal, with $x+z=2y$. For odd $m$ and $S\subseteq(\Zm)^d$ this means $x\ne z$, since $x=z$ forces $2y=2x$
and hence $y=x$.

\begin{proposition}\label{prop:relaxed}
Let $m$ be odd and $S\subseteq(\Zm)^d$. The following are equivalent.
\begin{enumerate}[label=\textup{(\arabic*)}]
\item $S$ is reduced to the empty set by the relaxed rule: repeatedly delete every point that is missing
from at least one of the three positions of the non-constant $3$-progressions in the current set.
\item $S$ is non-mid-point reducible.
\item There is a function $h\colon S\to\mathbb{R}$ with $h(x)+h(z)>2h(y)$ for every non-constant
$3$-progression $(x,y,z)$ in $S$.
\end{enumerate}
\end{proposition}

\begin{proof}
(2)$\Rightarrow$(1): A non-mid-point of the current set is missing from the middle position, so every
round of the relaxed rule deletes all non-mid-points, and at least one point as long as the current set,
which is reducible, is non-empty.

(1)$\Rightarrow$(2): Suppose that $S$ is not reducible. Then there is a non-empty set $C\subseteq S$ in which
every point is a mid-point. For every $y\in C$ choose a progression $(x_y,y,z_y)$ with distinct
$x_y,z_y\in C$, and form the directed graph on $C$ with the edges $y\to x_y$ and $y\to z_y$. Let $D$ be a
strongly connected component of this graph that no edge leaves. For $y\in D$ the points $x_y$ and $z_y$ lie
in $D$, so every point of $D$ is the middle term of a non-constant progression in $D$. Since $x_y\ne y$, the
component $D$ has at least two points, so every point $p$ of $D$ is the end of an edge that starts in $D$,
that is, $p\in\{x_y,z_y\}$ for some $y\in D$. As $(z_y,y,x_y)$ is a progression as well, $p$ occurs in the
first and in the third position of non-constant progressions in $D$. So as long as $D$ is contained in the
current set, the relaxed rule deletes no point of $D$, and it never empties $S$.

(2)$\Rightarrow$(3): Let $s_1,\dots,s_N$ be a removal order of $S$ and put $h(s_i)=3^{N-i}$. If $(x,y,z)$ is a
non-constant progression and $y=s_j$, then $y$ is a non-mid-point of $\{s_j,\dots,s_N\}$, so one of $x,z$ is
a point $s_i$ with $i<j$, and $h(x)+h(z)>3^{N-i}\ge3h(y)$.

(3)$\Rightarrow$(2): Let $S'$ be a non-empty subset of $S$ and $y$ a point of $S'$ where $h$ is largest. If
$2y=x+z$ with distinct $x,z\in S'$, then $(x,y,z)$ is a non-constant progression, so $h(x)+h(z)>2h(y)$,
against the choice of $y$. Hence $y$ is a non-mid-point of $S'$.
\end{proof}

The argument for (1)$\Rightarrow$(2) uses only that the end-points of a progression are distinct and differ
from its middle term, so it is not special to groups. The function $h$ in~(3) is a finite analogue of the
potentials of Section~\ref{sec:behrend}. A variant of it gives a short proof of
Proposition~\ref{prop:eps} for all moduli, odd and even.

\begin{remark}\label{rem:eps24}
Let $m\ge2$ and let $S\subseteq(\Zm)^d$ be non-mid-point reducible. Call two points $a,b$ of $S$ equivalent
if $2a=2b$. Whether a point $y$ is a mid-point of a set depends only on $2y$. So if a class contains a
non-mid-point of a subset of $S$ that contains the whole class, then all its points are non-mid-points of
this subset. Hence $S$ can be reduced by removing whole classes $C_1,\dots,C_N$, one after the other, such
that the points of $C_i$ are non-mid-points of $C_i\cup\dots\cup C_N$. Put $h(s)=3^{N-i}$ for $s\in C_i$. If
$x+z=2y$ with $x\ne z$ in $S$ and $y\in C_j$, then one of $x,z$ lies in a class $C_i$ with $i<j$, so
$h(x)+h(z)>2h(y)$. If $x+z=2y$ with $x=z$, then $x$ and $y$ are equivalent, so $h(x)+h(z)=2h(y)$.

Now let $x,y,z$ be $k$-tuples of elements of $S$ in which every element of $S$ occurs equally often, with
$x+z=2y$ and $x\ne z$. Then the sums $\sum_ih(x_i)$, $\sum_ih(y_i)$ and $\sum_ih(z_i)$ are equal. On the
other hand $h(x_i)+h(z_i)\ge2h(y_i)$ for every $i$, with strict inequality whenever $x_i\ne z_i$, which
happens for some $i$. Summing over $i$ gives a contradiction. So the set of these $k$-tuples is
progression-free. If $k$ is a multiple of $|S|$, it has at least $|S|^k/(k+1)^{|S|}$ elements, because the
tuples in which all elements occur equally often form the largest of the at most $(k+1)^{|S|}$ classes of
tuples with given numbers of occurrences. Padding with constant coordinates gives all dimensions $n$, and
this proves Proposition~\ref{prop:eps} with the exponent $n/d-O(\log n)$ for fixed $m$ and $d$. Unlike the
potentials of Section~\ref{sec:behrend}, the function $h$ is not bounded independently of $m$.
\end{remark}

Proposition~\ref{prop:relaxed} concerns progressions of length three; we make no claim for longer
progressions.

\section*{Acknowledgements}

This work was carried out with the help of Claude Fable 5, 5.2 and Opus 5.5. Claude contributed to the search for
the constructions, to the proofs, and to the computer programs used for the verification. ChatGPT was used for
proofreading and revisions, with further ideas and improvements. The first named author wishes to thank his
collaborators on previous papers on this subject. Both authors thank Michael H\"ausser for advice.
C.~Elsholtz was supported by the Austrian
Science Fund (FWF), project ArithRand I~4945.

\end{document}